\documentclass[final,leqno,onefignum,onetabnum]{siamltex1213}
\usepackage{cite,extarrows,bm,amssymb,amsmath,enumerate,color,txfonts,mathrsfs,tipa,multirow,epsfig,indentfirst}%
\usepackage{float}

\def\be{\begin{equation}}
\def\en{\end{equation}}
\def\beq{\begin{eqnarray}}
\def\eq{\end{eqnarray}}
\def\beqx{\begin{eqnarray*}}
\def\eqx{\end{eqnarray*}}

\makeatletter
  \newcommand\figcaption{\def\@captype{figure}\caption}
  \newcommand\tabcaption{\def\@captype{table}\caption}
\makeatletter
\newtheorem{remark}{Remark}[section]

\title{ Global space-time Trefftz DG schemes for the time-dependent linear wave equation 
\thanks{This work was completed during the author's visit to the University of Colorado at Boulder. The author was supported by the China Scholarship Council and Shandong Provincial Natural Science Foundation under the grant ZR2020MA046.}}
\author{LONG YUAN
\thanks{ College of Mathematics and Systems Science, Shandong
University of Science and Technology,
Qingdao 266590, China (yuanlong@lsec.cc.ac.cn).}
}

\begin{document}

\maketitle
\begin{abstract}
In this paper we are concerned with Trefftz discretizations of the time-dependent linear wave equation in anisotropic media in arbitrary space dimensional domains $\Omega \subset \mathbb{R}^d~ (d\in \mathbb{N})$. We propose two variants of the Trefftz DG method, define novel plane wave basis functions based on rigorous choices of scaling transformations and coordinate transformations,  and prove that the corresponding approximate solutions possess optimal-order error estimates with respect to the meshwidth $h$ and the condition number of the coefficient matrices, respectively.  Besides, we propose the global Trefftz DG method combined with local DG methods to solve the time-dependent linear nonhomogeneous wave equation in anisotropic media. In particular, the error analysis holds for the (nonhomogeneous) Dirichlet, Neumann, and mixed boundary conditions from the original PDEs. Furthermore, a strategy to discretize the model in heterogeneous media is proposed.
The numerical results verify the validity of the theoretical results, and show that the resulting approximate solutions possess high accuracy.
\end{abstract}


\pagestyle{myheadings} \thispagestyle{plain} \markboth{ LONG YUAN}{Trefftz DG methods for anisotropic wave equations}

\section{Introduction}
The idea at the heart of Trefftz method, which are named after the seminal work \cite{Tre} of E. Trefftz, is to choose the Trefftz approximation functions from a class of piecewise solutions of the same governing partial differential equation (PDE) without boundary conditions. Trefftz methods turned out to be particularly effective, and popular, for wave propagation problems in time-harmonic regime at medium and high frequencies, where the oscillatory nature of the solutions makes standard methods computationally too expensive; see the recent survey \cite{HMPsur} and references therein. The Trefftz method has an important advantage over Lagrange finite
elements for discretization of the Helmholtz equation and time-harmonic Maxwell equations \cite{ref21,pwdg,hy2,HMPsur,hmm,KMW,peng,peng2,yuan}: to achieve the same accuracy, relatively smaller degrees of freedom are enough in the plane wave-type methods owing to the particular choice of the basis functions that (may approximately) satisfy the considered PDE without boundary conditions.

Much work has been devoted to Trefftz discontinuous Galerkin (DG) methods for time-dependent linear isotropic wave phenomena, see in particular \cite{BMPS, egger,egger2, Kretzschmar2, Kretzschmar,Kretzschmar3, MP,PSSW}. A space-time Trefftz discontinuous Galerkin method for the first-order transient acoustic wave equations in arbitrary space dimensions is proposed and systematically studied in \cite{MP}. 
A Trefftz DG method for time-dependent electromagnetic problems has been analysed in \cite{Kretzschmar,Kretzschmar3} in one space dimension, and then it has been extended to three-dimensional time dependent Maxwell's equations in \cite{egger,egger2,Kretzschmar2}. Besides, the recent work on explicit marching-type space-time schemes, which are so-called {\it tent-pitched space-time meshes}, are fully investigated by \cite{falk,glpal,glpal2,PSSW}, where the PDE are explicitly evolved from the ``bottom" to the ``top" of the space-time cylinder element by element.

Recently, a class of space-time DG discretizations of the linear isotropic acoustic wave equation in two space dimensions in polygonal domains occupied by possibly heterogeneous media are addressed in \cite{BMPS}.  The consistency analysis of the space-time discretization are generalized to non-Trefftz discrete spaces, and the realistic setting of solutions exhibiting spatial point singularities is allowed, where high $h-$convergence rates of the approximations generated by the space-time scheme with local corner mesh refinement on the spatial domain still hold. Moreover, the consistency error bounds in mesh-dependent norm holds true without any time-step size constraint.

Anisotropy can result from periodic layering of fine layers \cite{Carcione}, preferential alignment of fractures and cracks \cite{Sayers}. Anisotropy may greatly influence seismic wave propagation, seismic data acquisition and subsequent data analysis and processing procedures \cite{Tsvankin}. It is therefore important to design accurate and efficient numerical methods for modeling wave propagation in anisotropic media. An acoustic wave equation for anisotropic media in \cite{Alkhalifah} is introduced to describe a wave type that propagates at speeds slower than the P-wave for a positive anisotropy parameter. An improved rotated staggered-grid finite-difference method in \cite{gao} with fourth-order temporal accuracy has been developed to solve elastic-wave modeling in anisotropic media, where the symmetry axes of anisotropy are not aligned with the coordinate axes. A compensated-amplitude vertical transverse isotropic least-squares reverse time migration method in \cite{qu} is adopted to correct the anisotropy effect and compensate amplitude attenuation. \cite{Zhan} develops a closed expression of Riemann solvers for the discontinuous Galerkin time domain method, applied to wave propagation modeling in distinct anisotropic material properties. Recently, the PWDG methods \cite{yuan} have been developed to solve Helmholtz equation and time-harmonic Maxwell equations in three-dimensional anisotropic media.

In this paper we construct a global space-time Trefftz DG scheme for the linear {\it anisotropic} wave equation in arbitrary space dimensional domains $\Omega \subset \mathbb{R}^d~ (d\in \mathbb{N})$. In order to build better convergence results, we have to carefully define plane wave basis functions by rigorous choices of the scaling transformations and the coordinate transformations. We prove that the approximate solutions generated by the proposed method possess satisfactory and optimal-order error estimates with respect to meshwidth $h$ and the condition number $\rho$ of the coefficient matrices, respectively. Besides, we propose an alternative standard Trefftz DG method with almost the same computational cost and with the same convergence order with respect to $h$ and $\rho$ in Section 6.

Numerical results indicates that, the consistency error bounds in the mesh-dependent  $|||\cdot|||_{\text{DG}}-$norm and mesh-independent $L^2(\Omega\times\{T\})-$norm for Trefftz DG methods are optimal with respect to $h$ and $\rho$, respectively. Moreover, numerical experiments in Section 8.3 show that the approximations generated by the Trefftz DG methods are clearly more accurate than that generated by the high-order DG finite element method.

Since Trefftz basis functions on each element are solutions of the {\it homogeneous} wave equations  without boundary conditions, the Trefftz methods can not be directly applied to  discretizations of the {\it nonhomogeneous} wave equations. Motivated by the coupled discontinuous Galerkin formulation developed in \cite{hy3}, we develop the global Trefftz DG method combined with overlapping local DG method. Numerical results indicates that, the consistency error bounds in the mesh-dependent  $|||\cdot|||_{\text{DG}}-$norm and mesh-independent $L^2(\Omega\times\{T\})-$norm for Trefftz DG method combined with {\it overlapping} local DG are optimal with respect to $h$ and $\rho$, respectively, in the presence of the nonhomogeneous source and the anisotropic media. Besides, we propose another alternative to define nonhomogeneous local problems in each {\it nonoverlapping} time slab, and solve it by the space-time DG method. The resulting residue problem on the global solution domain is still solved by the Trefftz DG method. We call the new method as `` Trefftz DG method combined with nonoverlapping local DG".  Numerical results in Section 8.4 indicates that, Trefftz DG method combined with nonoverlapping local DG is comparable to  Trefftz DG method combined with overlapping local DG.

Furthermore,  in order to make our proposed method adaptive to the model in heterogeneous media where $A$ is a piecewise-constant positive definite matrix, a second-best strategy that the computational space domain $\Omega$ is directly partitioned is employed such that the mesh ${\cal T}_{h_{{\bf x}}}^{\bf x}=\{K_{\bf x}\}$ satisfies the shape regular and quasi uniform conditions. The error estimates of corresponding Trefftz discontinuous Galerkin approximations are proved in section \ref{piecesec} and numerical results are reported in section \ref{heterosec}.

Comparing against the space-time DG method of \cite{BMPS} for isotropic wave equations which employs the piecewise-polynomial discrete space, our global (resp. local) discrete space is available for the Trefftz (resp. DG finite element) space in arbitrary space dimensions, and the Dirichlet, Neumann, and mixed boundary conditions on space-time domain boundary are considered. In particular, there is no constraint on the data of homogeneous Neumann boundary conditions from the original
initial boundary value problem (IBVP), thus the analysis holds for the nonhomogeneous Neumann boundary conditions from the original IBVP, which has been also verified by the numerical tests in Section 8.3.

The analysis framework presented in this paper is borrowed from \cite{MP}. The extending steps to the anisotropic case consist of establishing the variational formulation in section 3,  constructing anisotropic Trefftz basis function spaces in section 4, identifying mesh skeleton norms in section 5 on the Trefftz function space for which the bilinear form defining the method is coercive, which allows us to prove well-posedness and error estimates in these norms. Moerover, in comparison to most existing Trefftz methods \cite{Kretzschmar2,Kretzschmar,MP} for the isotropic wave equations, the proposed variational formulation with three relaxation parameters is applied to the anisotropic model with Dirichlet, Neumann, and mixed boundary conditions; combined with the local DG method on auxiliary smooth subdomains, the new method can generate the approximations with the spectral convergence orders for the nonhomogeneous case; the space-time domain partition ${\cal T}_h$ is obtained as the tensor product of space and time mesh grids ${\cal T}_{h_{{\bf x}}}^{\bf x}$ and ${\cal T}_{h_t}^t$, the space-time mesh $\hat{\cal T}_{\hat h}$ of the transformed space-time domain $\hat Q$ satisfies quasi-uniform assumption, and the space mesh grid $\hat{\cal T}_{\hat h_{\hat{\bf x}}}^{\hat{\bf x}}$ of $\hat\Omega$ is shape regular and quasi-uniform.

The paper is organized as follows: In Section 2, we state the initial boundary value problem for the acoustic wave equation in both first- and second-order formulation. Section 3 describes the proposed method for the homogeneous PDEs.
In Section 4, we explain how to discretize the resulting variational problems. Section 5 provides the desired error estimates for the approximate solutions.
In Section 6, we propose an alternative standard Trefftz DG method. In order to solve the nonhomogeneous and anisotropic model, we develop a global Trefftz DG method combined with overlapping local DG method for the nonhomogeneous in Section 7. In Section 8, we introduce another strategy to discretize the model in heterogeneous media. Finally, we report some numerical results to confirm the effectiveness of the
proposed method .

\section{Considered model}
 We consider the first order acoustic wave IBVP posed on
 a space-time domain $Q=\Omega\times I$, where $\Omega \subset \mathbb{R}^d~ (d\in \mathbb{N})$ is an open bounded Lipschitz polytope and $I=(0,T),T>0$.
 ${\bf n}^x_{\Omega}$ is an outward-pointing unit normal vector on $\partial\Omega$.
The boundary of $\Omega$ denoted by $\Gamma$, is divided in two parts, with mutually disjoint interiors, denoted $\Gamma_D$ or $\Gamma_N$ corresponding to Dirichlet and Neumann boundary conditions, respectively. The model reads as
  \begin{equation} \label{model}
\left\{ \begin{aligned}
     &   A^{\frac{1}{2}}\nabla v + \frac{\partial {\bm \sigma}}{\partial t}={\bf 0} & \text{in} \quad Q, \\
      & \nabla\cdot A^{\frac{1}{2}}{\bm \sigma} + c^{-2} \frac{\partial v}{\partial t}=0& \text{in} \quad Q, \\
      &  v(\cdot,0)=v_0, \quad {\bm \sigma}(\cdot,0)= {\bm \sigma}_0 & \text{on} \quad \Omega, \\
       &    v = g_D & \text{on} \quad \Gamma_D\times[0,T],
       \\
      &   A^{\frac{1}{2}}{\bm \sigma} \cdot {\bf n}^x_{\Omega} = g_N & \text{on} \quad \Gamma_N\times[0,T].
                          \end{aligned} \right.
                          \end{equation}
Here $v_0, {\bm \sigma}_0, g_D, g_N$ are the given source data, $c$ is the wave speed, which is constant in the whole space domain $\Omega$ and independent of time $t$. $A$ is positive definite matrix independent of ${\bf x}$ and $t$ (See section \ref{piecesec} for the case of piecewise constant matrices). The gradient $\nabla$ and divergence $\nabla\cdot$ operators are meant in the space variable ${\bf x}$ only.

If there exists a scalar field $U_0$ such that ${\bm \sigma}_0 =-A^{\frac{1}{2}}\nabla U_0$, then IBVP (\ref{model}) is equivalent to the following second order scalar wave equation, by setting $v = \frac{\partial U}{\partial t}$ and ${\bm \sigma}=-A^{\frac{1}{2}}\nabla U$,
 \begin{equation} \label{scalarmodel}
\left\{ \begin{aligned}
       & -\nabla\cdot (A\nabla U)  + c^{-2} \frac{\partial^2 U}{\partial t^2}=0 & \text{in} \quad Q, \\
      &  \frac{\partial U}{\partial t}(\cdot,0)=v_0, \quad U(\cdot,0)= U_0 & \text{on} \quad \Omega, \\
            &     \frac{\partial U}{\partial t}  = g_D & \text{on} \quad \Gamma_D\times[0,T],
       \\
      &   -A\nabla U \cdot {\bf n}^x_{\Omega} = g_N & \text{on} \quad \Gamma_N\times[0,T].
                          \end{aligned} \right.
                          \end{equation}

Let the time domain $(0,T)$ be divided into $N\in \mathbb{N}$ intervals $I_n (1\leq n \leq N)$ composing a partition ${\cal T}_{h_t}^t$, with
$$I_n=(t_{n-1},t_n),~~h_n=t_n-t_{n-1}=|I_n|, ~~h_t=\mathop\text{max}\limits_{1\leq n \leq N}h_n.$$
Let us introduce the following notation for the time slabs and the partial cylinders, respectively,
$$D_n = \Omega\times I_n, ~~Q_n = \Omega\times(0,t_n), ~~1\leq n \leq N.$$

For each $1\leq n \leq N$, we introduce a same polygonal finite element mesh ${\cal T}_{h_{{\bf x}}}^{\bf x}=\{K_{\bf x}\}$ of the spatial domain $\Omega$ with
$$h_{K_{\bf x}}=\text{diam}K_{\bf x}, ~~h_{{\bf x}}=\mathop\text{max}\limits_{K_{\bf x}\in {\cal T}_{h_{\bf x}}^x}h_{K_{\bf x}}.$$

Then the space-time domain $Q=\Omega\times(0,T)$ can be partitioned with a finite element mesh ${\cal T}_h$ given by
$${\cal T}_h=\{K=K_{\bf x}\times I_n, K_{\bf x}\in {\cal T}_{h_{{\bf x}}}^{\bf x}, 1\leq n \leq N\}.$$
Here ${\cal T}_h$ is a tensor product mesh.
Besides, we define the time-truncated mesh
$${\cal T}_h(Q_n) = \{K\in {\cal T}_h, K\subset Q_n \}, ~1\leq n \leq N$$
and
$${\cal T}_h(D_n) = \{K\in {\cal T}_h, K\subset D_n \}, ~1\leq n \leq N.$$

Assume the space-time grid \( {\cal T}_h\) satisfies the assumptions presented in Section 4 of \cite{BMPS}: on an internal face $F=\partial K_1\bigcap \partial K_2$, either
\begin{equation} \label{spacetimemesh}
\left\{ \begin{aligned}
     &   {\bf n}_F^{\bf x}=0  & \text{and} ~F ~\text{is called ``space-like" face, or} \\
      &  {\bf n}_F^t=0  & \text{and} ~F ~\text{is called ``time-like" face}, \\
                          \end{aligned} \right.
                          \end{equation}
where $({\bf n}_F^{\bf x},{\bf n}_F^t)$ is a unit vector of the face $F$. On space-like faces, by convention, we choose ${\bf n}_F^t>0$, which means that the unit normal vector $({\bf n}_F^{\bf x},{\bf n}_F^t)$ points towards future time. Moreover, all time-like faces are of the form $F=F_{\bf x}\times F_t$ with $h_{F_{\bf x}}=|F_{\bf x}|$ and $h_{F_t}=|F_t|$; we recall that $F_t=I_n, 1\leq n \leq N$. Finally, We denote the outward-pointing unit
normal vector on $\partial K$ by $({\bf n}_K^{\bf x},{ n}_K^t)$ .

We denote by $\mathcal{F}_h = \bigcup\limits_{K\in{\cal T}_h}\partial K$ the
skeleton of the mesh,
 by $\mathcal{F}_h^{\text{space}}$ the union of the internal space-like faces,  and by $\mathcal{F}_h^{\text{time}}$ the union of the internal time-like faces, respectively.  Set $\mathcal{F}_h^0=\Omega\times \{t=0\}$, $\mathcal{F}_h^T=\Omega\times \{t=T\}$, $\mathcal{F}_h^D=\Gamma_D\times [0,T]$ and $\mathcal{F}_h^N=\Gamma_N\times [0,T]$.

Let $w$, ${\bm \tau}$ and $M$ be a piecewise smooth function, vector field and matrix function on ${\cal T}_h$, respectively.
On $F=\partial K_1\bigcap\partial K_2$, we define
\begin{eqnarray}
& \text{the averages:} ~~\{\{ w \}\}:= \frac{w_{|K_1}+w_{|K_2}}{2}, ~~\{\{
{\bm \tau}\}\} := \frac{{\bm \tau}_{|K_1} + {\bm \tau}_{|K_2}}{2},
\cr &
\text{space normal jumps:} ~\llbracket {w} \rrbracket_{\bf N} :=
w_{|K_1}{\bf n}^{\bf x}_{K_1} + w_{|K_2}{\bf n}^{\bf x}_{K_2}, ~\llbracket {\bm \tau} \rrbracket_{\bf N} = {\bm \tau}_{|K_1}\cdot {\bf n}^{\bf x}_{K_1} + {\bm \tau}_{|K_2}\cdot {\bf n}^{\bf x}_{K_2},
\cr &
~\llbracket M \rrbracket_{\bf N} :=
M_{|K_1}{\bf n}^{\bf x}_{K_1} + M_{|K_2}{\bf n}^{\bf x}_{K_2},
\cr &
\text{time full jumps:} ~\llbracket w\rrbracket_t:= w_{|K_1}{ n}^t_{K_1} + w_{|K_2}{ n}^t_{K_2} = (w^- - w^+) { n}^t_F,
\cr &
\text{time full jumps:} ~\llbracket {\bm \tau} \rrbracket_t:= {\bm \tau}_{|K_1}{ n}^t_{K_1} + {\bm \tau}_{|K_2}{ n}^t_{K_2}=({\bm \tau}^- - {\bm \tau}^+){ n}^t_F.
\end{eqnarray}
 Here $w^-$ and $w^+$ denote the traces of the function $w$ from the adjacent elements at lower and higher times, respectively, and similarly for ${\bm \tau}^{\pm}$.

\section{The variational formulation}
Set local Trefftz space:
\beq
{\bf T}(K) &=& \bigg\{  (w,{\bm \tau}) \in H^1(K)^{1+d} s.t.  ~{\bm \tau}|_{\partial K}\in L^2(\partial K)^d, ~\frac{\partial w}{\partial t}, \nabla\cdot {\bm \tau} \in L^2(K),
 \cr &&\frac{\partial\bm \tau}{\partial t}, \nabla w\in L^2(K)^{d},
~ A^{\frac{1}{2}}\nabla w + \frac{\partial {\bm \tau}}{\partial t}={\bf 0} ,
  \nabla \cdot A^{\frac{1}{2}}{\bm \tau} + c^{-2} \frac{\partial w}{\partial t}=0
  \bigg\} ~\quad \forall K\in {\cal T}_h,
\eq
and set global Trefftz space:
\beq
{\bf T}({\cal T}_h) = \bigg\{  (w,{\bm \tau}) \in L^2(Q)^{1+d} ~s.t.  ~(w|_K,{\bm \tau}|_K)  \in {\bf T}(K) ~\forall K\in {\cal T}_h \bigg\}.
\eq

To derive the Trefftz-DG  variational formulation, we multiply the first two equation of (\ref{model}) with test fields ${\bm \tau}$ and $w$ and integrated by parts on each $K\in {\cal T}_h$:
\beq \label{stavaria0}
 -\int_K && \bigg( v (\nabla\cdot A^{\frac{1}{2}} {\bm \tau} + c^{-2}\frac{\partial w}{\partial t}) + {\bm \sigma} \cdot(A^{\frac{1}{2}}\nabla w + \frac{\partial \bm\tau}{\partial t})\bigg)dV
 \cr + && \int_{\partial K}\bigg( (vA^{\frac{1}{2}} {\bm \tau} + w A^{\frac{1}{2}}{\bm \sigma} )\cdot {\bf n}^x_K + ( {\bm \sigma} \cdot {\bm \tau} +   c^{-2} v w )n^t_K \bigg)dS = 0, ~ \forall (w,{\bm \tau}) \in {\bf T}({\cal T}_h).
\eq
Replacing the traces of $v$ and ${\bm\sigma}$ on the mesh skeleton by the single-valued numerical fluxes $\check{v} $ and $\check{\bm\sigma}$, we have
\be \label{stavaria}
  \int_{\partial K}\bigg( \check{v}( A^{\frac{1}{2}} {\bm \tau}\cdot {\bf n}^x_K + c^{-2} w n^t_K)  + \check{\bm \sigma} \cdot ( w A^{\frac{1}{2}} {\bf n}^x_K + {\bm \tau} n^t_K) \bigg)dS = 0.
\en
Define the numerical fluxes as follows.

\begin{equation} \label{standflux} \nonumber
\check{v}=\left\{ \begin{aligned}
     &  v^- \\
     & v \\
     & v_0 \\
     & \{\{v\}\} + \beta \llbracket  A^{\frac{1}{2}}{\bm\sigma} \rrbracket_N \\
     & g_D \\
     & v + \beta (A^{\frac{1}{2}}{\bm \sigma} \cdot {\bf n}^x_{\Omega} - g_N)
  \end{aligned} \right.
  \check{\bm \sigma}= \left\{ \begin{aligned}
     &  {\bm \sigma}^- &\quad {\text on} ~ {\mathcal{F}}_{ h}^{\text{space}},\\
     & {\bm \sigma} &\quad {\text on} ~ {\mathcal{F}}_{ h}^{T}, \\
     & {\bm \sigma}_0 &\quad {\text on} ~ {\mathcal{F}}_{ h}^{0}, \\
     &  \{\{\bm \sigma\}\} + \alpha\llbracket A^{\delta}v \rrbracket_N &\quad {\text on} ~ {\mathcal{F}}_{ h}^{\text{time}}, \\
     & \bm \sigma + \alpha(v-g_D) A^{\delta} {\bf n}^x_{\Omega} &\quad {\text on} ~ {\mathcal{F}}_{ h}^D, \\
     & g_N A^{-\frac{1}{2}} {\bf n}^x_{\Omega} &\quad {\text on} ~ {\mathcal{F}}_{ h}^N.
  \end{aligned} \right.
                          \end{equation}
The stabilization parameters $\alpha \in L^{\infty}({\mathcal{F}}_{ h}^{\text{time}}\bigcup {\mathcal{F}}_{ h}^D)$, $\beta \in L^{\infty}({\mathcal{F}}_{ h}^{\text{time}}\bigcup {\mathcal{F}}_{ h}^N)$ are positive constant on each time-like face. $\delta \in L^{\infty}({\mathcal{F}}_{ h}^{\text{time}}\bigcup {\mathcal{F}}_{ h}^D)$ is constant on each time-like face, and its best choice will be given in Lemma \ref{helmimportl}.

By summing the elemental DG equation over the element $K\in {\bf T}({\cal T}_h)$ and using the defined fluxes,  we can obtain the Trefftz-DG  variational formulation:
Find $(v,{\bm \sigma}) \in {\bf T}({\cal T}_h)$ such that
\be \label{trefftzcontivaria}
\mathcal{A}(v,{\bm \sigma};w,{\bm \tau}) =
\ell(w,{\bm \tau}) \quad \forall (w,{\bm \tau}) \in {\bf T}({\cal T}_h),
\en
where
\beq \label{adg1}
\mathcal{A}(v,{\bm \sigma};w,{\bm \tau})&= & \int_{\mathcal{F}_h^{\text{space}}}\big(c^{-2}v^-\llbracket w\rrbracket_t + {\bm \sigma}^-\cdot \llbracket {\bm \tau} \rrbracket_t  \big)~d{\bf x}
+  \int_{\mathcal{F}_h^{\text{time}}}
\bigg(\{\{v\}\} \llbracket {A^{\frac{1}{2}}\bm \tau} \rrbracket_{\bf N} +
 \{\{ {\bm \sigma} \}\} \cdot \llbracket A^{\frac{1}{2}}w \rrbracket_{\bf N}
 \cr &+&  \alpha\llbracket A^{\delta}v \rrbracket_{\bf N} \cdot \llbracket A^{\frac{1}{2}}w \rrbracket_{\bf N}
  + \beta \llbracket A^{\frac{1}{2}}{\bm \sigma} \rrbracket_{\bf N} \llbracket A^{\frac{1}{2}}{\bm \tau} \rrbracket_{\bf N}  \bigg)~dS
    \cr & + &  \int_{\mathcal{F}_h^T} \big( c^{-2}v w + {\bm \sigma} \cdot {\bm \tau}\big)~d{\bf x}
    +  \int_{\mathcal{F}_{h}^D} ( {\bm \sigma} \cdot w A^{\frac{1}{2}}{\bf n}^x_{\Omega} + \alpha v w A^{\delta}{\bf n}^x_{\Omega}\cdot A^{\frac{1}{2}}{\bf n}^x_{\Omega}) ~dS
    \cr & + &   \int_{\mathcal{F}_{h}^N} \bigg(  v (A^{\frac{1}{2}}{\bm \tau}\cdot {\bf n}^x_{\Omega}) + \beta  (A^{\frac{1}{2}}{\bm \sigma} \cdot {\bf n}^x_{\Omega})
  (A^{\frac{1}{2}}{\bm \tau}\cdot {\bf n}^x_{\Omega}) \bigg)~dS,
 \eq
 and
\beq \label{adglf1}
 \ell(w,{\bm \tau}) &=& \int_{\mathcal{F}_h^0}\big( c^{-2}v_0 w + {\bm \sigma}_0 \cdot {\bm \tau}\big)~d{\bf x}
   +   \int_{\mathcal{F}_h^N} g_N \bigg( \beta A^{\frac{1}{2}}{\bm \tau}\cdot {\bf n}^x_{\Omega}  - w \bigg) ~ dS
    \cr & + &  \int_{\mathcal{F}_h^D} \alpha g_D w A^{\delta}{\bf n}^x_{\Omega}\cdot A^{\frac{1}{2}}{\bf n}^x_{\Omega} ~ dS
 - \int_{\mathcal{F}_h^D} g_D A^{\frac{1}{2}}{\bm \tau}\cdot {\bf n}^x_{\Omega} ~ dS.
 \eq

\begin{remark}
The choice of numerical fluxes $(\check{v}, \check{\bm \sigma})$
defined on ``time-like" interfaces and boundary faces is such that
the Trefftz DG formulation is consistent; namely, if $(v,
\bm\sigma)\in H^1(Q)$ solves (\ref{model}), then it satisfies (\ref{trefftzcontivaria}). In particular, the flux $\check{\bm \sigma}$ satisfies $A^{\frac{1}{2}}\check{\bm \sigma} \cdot {\bf n}^x_{\Omega} = g_N$ on ${\mathcal{F}}_{ h}^N$ coinciding with the Neumann boundary condition satisfied by the exact solution $(v,\bm\sigma)$ of the IBVP (\ref{model}).
\end{remark}

\section{Discretization of the variational problems}
The proposed Trefftz DG method for (\ref{model}) depends on two transformations.

\subsection{A coordinate transformation and a scaled transformation}
 Since $A$ is positive definite matrix, there exists an orthogonal matrix $P$ and a diagonal positive definite matrix $\Lambda = \text{diag} (\lambda_1,\lambda_2,\cdots,\lambda_d)$ such that $A= P^T \Lambda P$, where $\lambda_i\leq\lambda_{i+1} (1\leq i\leq d-1)$ and the superscript $T$ denotes matrix transposition. Set $\lambda_{\text{min}}=\lambda_1$ and $\lambda_{\text{max}}=\lambda_d$.
Of course, we can assume that $\text{det}(P)=1$. Define a coordinate transformation:
\be\label{tran1}
 \hat{\bf x} = \Lambda^{-\frac{1}{2}}P{\bf x} \xlongequal{\Delta} S{\bf x}, \quad S=\Lambda^{-\frac{1}{2}}P.
\en

Under the coordinate transformation (\ref{tran1}), let $\hat{\Omega}$ and $\hat Q$ denote the images of $\Omega$ and $Q$, respectively, and denote by $\hat{\cal T}_{{\hat h}_{\hat{\bf x}}}^{\hat{\bf x}}=\{\hat K_{\hat{\bf x}}\}$ the transformed finite element mesh of the spatial domain $\hat\Omega$ with
$$\hat h_{\hat K_{\hat{\bf x}}}=\text{diam}\hat K_{\hat {\bf x}}, ~~\hat h_{\hat{\bf x}}=\mathop\text{max}\limits_{\hat K_{\hat{\bf x}}\in \hat{\cal T}_{{\hat h}_{\hat {\bf x}}}^{\hat x}} \hat h_{\hat K_{\hat {\bf x}}}.$$
Furthermore, the transformed space-time domain $\hat Q=\hat\Omega\times(0,T)$ can be partitioned with a finite element mesh $\hat{\cal T}_{\hat h}$ given by
$$\hat{\cal T}_{\hat h}=\{\hat K=\hat K_{\hat{\bf x}}\times I_n, \hat K_{\hat{\bf x}}\in \hat{\cal T}_{\hat h_{\hat{\bf x}}}^{\hat{\bf x}}, 1\leq n \leq N\}.$$

Assume that the space-time mesh satisfies the condition
\be
\hat h = \mathop\text{max}\big\{ \hat h_{\hat{\bf x}}, \mathop\text{max}\limits_{n^{'}=1,\cdots,n} ch_{n^{'}} \big\} \leq \hat\rho \mathop\text{min}\big\{ \hat h_{\hat{\bf x}}, \mathop\text{min}\limits_{n^{'}=1,\cdots,n} c h_{n^{'}} \big\},
\en
for each discrete time $t_n$ and some $\hat\rho>1$.

We denote by $\hat{\mathcal{F}}_{\hat h} = \bigcup\limits_{\hat K\in {\hat{\cal T}}_{\hat h}} \partial \hat K$ the
skeleton of the mesh,
 by $\hat{\mathcal{F}}_{\hat h}^{\text{space}}$ the union of the internal space-like faces,  and by $\hat{\mathcal{F}}_{\hat h}^{\text{time}}$ the union of the internal time-like faces, respectively.
 We use $\hat{\bf n}$ to denote the unit outer normal vector on the boundary of each element \(\hat K_{\hat{\bf x}}\). Denote by $\hat\Gamma_D$ and $\hat\Gamma_N$ the images of $\Gamma_D$ and $\Gamma_N$ under the coordinate transformation (\ref{tran1}), respectively.
Set $\hat{\mathcal{F}}_{\hat h}^0=\hat\Omega\times \{t=0\}$,
$\hat{\mathcal{F}}_{\hat h}^T=\hat\Omega\times \{t=T\}$,
$\hat{\mathcal{F}}_{\hat h}^D=\hat\Gamma_D\times [0,T]$ and $\hat{\mathcal{F}}_{\hat h}^N=\hat\Gamma_N\times [0,T]$.

Denote by $\nabla_h$ and $\hat\nabla_{\hat h}$ the element application of the spacial gradient operator $\nabla=(\frac{\partial}{\partial x_1} ~\frac{\partial}{\partial x_2} ~\cdots ~\frac{\partial}{\partial x_d} )^T$ and $\hat\nabla=(\frac{\partial}{\partial \hat x_1} ~\frac{\partial}{\partial \hat x_2} ~\cdots~\frac{\partial}{\partial \hat x_d} )^T$, respectively.
Define the spacial Laplace operator $\hat\triangle$ on $\hat{\Omega}$ by $\hat\triangle=\frac{\partial^2}{\partial\hat x_1^2}+\frac{\partial^2}{\partial\hat x_2^2}+\cdots +\frac{\partial^2}{\partial\hat x_d^2}$.

Define the scaled fields $(\hat v,\hat{\bm \sigma})$ as
 \be \label{firstorderrela}
\hat v = v, ~~ \hat{\bm \sigma}= P {\bm \sigma}.
\en

By some patient calculation in ``Appendix",
we can obtain the following relationships transforming ``anisotropic" into ``isotropic":
\be \label{initanisotoiso}
A^{\frac{1}{2}}\nabla v = P^T \hat\nabla\hat v, ~~\text{and} ~~\nabla\cdot (A^{\frac{1}{2}}{\bm \sigma})=\hat\nabla \cdot \hat{\bm \sigma}.
\en
 Thus the anisotropic wave equation (\ref{model}) is transformed into the isotropic wave equation:
\begin{equation} \label{transmodel}
\left\{ \begin{aligned}
     &   \hat\nabla \hat v + \frac{\partial \hat{\bm \sigma}}{\partial t}={\bf 0} & \text{in} \quad \hat Q, \\
      & \hat\nabla\cdot \hat{\bm \sigma} + c^{-2} \frac{\partial \hat v}{\partial t}=0 & \text{in} \quad \hat Q. \\
                          \end{aligned} \right.
                          \end{equation}

 Conversely, if $(\hat v, \hat{\bm \sigma})$ satisfies the isotropic isotropic wave equation (\ref{transmodel}),  $(v,{\bm\sigma})$ defined by the inverse scaled transformation of (\ref{firstorderrela}) and the coordinate transformation (\ref{tran1}) satisfies the original anisotropic wave equation (\ref{model}).

\subsection{Anisotropic Trefftz basis function spaces}
In order to derive a finite dimensional Trefftz space ${\bf V}_h({\cal T}_h) \subset {\bf T}({\cal T}_h) $ satisfying the original anisotropic wave equation (\ref{model}),  we first give the definition of a discretized Trefftz space $\hat{\bf V}_{\hat h}(\hat{\cal T}_{\hat h})$ satisfying  isotropic wave equation (\ref{transmodel}). We refer the reader to \cite[Remark 13]{MP} for a detailed construction of $\hat{\bf V}_{\hat h}(\hat{\cal T}_{\hat h})$ by evolving in time polynomial initial conditions.

Assuming that the first order problem (\ref{model}) is derived from the second order problem (\ref{scalarmodel}), and define the polynomial Trefftz space for the second order problem from (\ref{transmodel}):
\be \label{secondtrefftz}
\hat{\mathbb{U}}^p(\hat K)=\{ \hat U \in \mathbb{P}^p(\hat K) ~ s.t. ~ -\Delta \hat U  + c^{-2} \frac{\partial^2 \hat U}{\partial t^2}=0 \},
\en
 where the subscript $p$ is related to the dimension of the local spaces.
 Denote some multi-index notation for ${\bm \alpha} \in \mathbb{N}_0^d$ by $|{\bm \alpha}|=\alpha_1+\cdots+\alpha_d$, $D^{\bm \alpha}\phi = \frac{\partial^{|{\bm \alpha}|} \phi}{\partial \hat x_1^{\alpha_1}\cdots \partial \hat x_d^{\alpha_d}}$,
${\bf x}^{\bm \alpha}=\hat x_1^{\alpha_1}\cdots \hat x_d^{\alpha_d}$, and for a space-time field $\phi$, by $D^{\alpha_t, \bm \alpha}\phi =
\frac{\partial^{\alpha_t + |{\bm \alpha}|} \phi}{\partial t^{\alpha_t}\partial \hat x_1^{\alpha_1}\cdots \partial \hat x_d^{\alpha_d}}$.
If the polynomial
\be \nonumber
\hat U(\hat{\bf x},t) = \sum\limits_{k\in \mathbb{N}_0,{\bm \alpha}\in \mathbb{N}_0^d, k+|{\bm \alpha}|\leq p } a_{k,{\bm \alpha}} t^k\hat{\bf x}^{\bm \alpha}
\en
with $a_{k,{\bm \alpha}}\in \mathbb{R}$ satisfies the second order wave equation in (\ref{secondtrefftz}), then the coefficients $a_{k,{\bm \alpha}}$ satisfy the recurrence
\be \nonumber
a_{k,{\bm \alpha}} = \frac{c^2}{k(k-1)}\sum_{m=1}^d (\alpha_m+2)(\alpha_m+1) a_{k-2,{\bm \alpha}+2{\bf e}_m},
\en
where ${\bf e}_m$ is the $m$th row of the identity matrix of order $d$.

In order to start the recursion, one can start by choosing polynomial basis functions $\{ \tilde{b}_1,\cdots, \tilde{b}_{C_{p+d}^d} \}$ for the space $\mathbb{P}^p(\mathbb{R}^d)$ for $k=0$ and $\{ \tilde{\tilde{b}}_1,\cdots, \tilde{\tilde{b}}_{C_{p-1+d}^d} \}$ for the space $\mathbb{P}^{p-1}(\mathbb{R}^d)$ for $k=1$. Then  a basis for $\hat{\mathbb{U}}^p(\hat K)$ can be defined such that either $U(\cdot,0)=\tilde b_{j}$ and $\frac{\partial U}{\partial t}(\cdot,0)=0$, or $U(\cdot,0)=0$ and $\frac{\partial U}{\partial t}(\cdot,0)=\tilde{\tilde{b}}_{j}$ for some $j$. It leads to the dimension of $\hat{\mathbb{U}}^p(\hat K)$ as
\be \nonumber
\text{dim} \hat{\mathbb{U}}^p(\hat K) = C_{p+d}^d + C_{p-1+d}^d,
\en
where $C_k^j=\frac{k!}{j!(k-j)!}$ for $j\leq k\in \mathbb{N}_0$.

We denote by $\hat b_j ( 1\leq j\leq \text{dim} \hat{\mathbb{U}}^p(K) )$ the basis functions of the space $\hat{\mathbb{U}}^p(\hat K)$.
Then, a Trefftz space $\hat{\mathbb{W}}^p(\hat K)$ for the first order system can be derived from 
\be \nonumber
\hat{\mathbb{W}}^p(\hat K) = \text{span}\bigg\{ (\frac{\partial \hat b_j}{\partial t},-\hat\nabla \hat b_j), ~\hat b_j\in  \hat{\mathbb{U}}^{p+1}(\hat K),~ 1\leq j\leq \text{dim} \hat{\mathbb{U}}^{p+1}(\hat K) \bigg\}.
\en
Since the constants in $\hat{\mathbb{U}}^{p+1}(\hat K)$ have no contribution to $\hat{\mathbb{W}}^p(\hat K)$,  the dimension of $\hat{\mathbb{W}}^p(\hat K)$ equals $\text{dim} \hat{\mathbb{U}}^{p+1}(\hat K)-1$.
Meanwhile we have the isotropic Trefftz space defined on $\hat{\cal T}_{\hat h}$
\be \nonumber
\hat{\bf V}_{\hat h}(\hat{\cal T}_{\hat h}) = \prod\limits_{\hat K\in \hat{\cal T}_{\hat h}} \hat{\mathbb{W}}^{p}(\hat K).
\en

By the coordinate transformation (\ref{tran1}) and (\ref{firstorderrela}), we get the anisotropic Trefftz space ${\mathbb{W}}^p(K)$ on $K\in {\cal T}_h$ for the first two equations of first order system (\ref{model})
\be
{\mathbb{W}}^p( K) = \text{span}\bigg\{ (\frac{\partial \hat b_j}{\partial t},-P^T\hat\nabla \hat b_j), ~\hat b_j\in  \hat{\mathbb{U}}^{p+1}(\hat K),~K=S^{-1}\hat K,~ 1\leq j\leq \text{dim} \hat{\mathbb{U}}^{p+1}(\hat K) \bigg\}.
\en
Furthermore,  we get the anisotropic Trefftz space defined on ${\cal T}_{h}$
\be
{\bf V}_{ h}({\cal T}_{ h}) = \prod\limits_{ K\in {\cal T}_{ h}} {\mathbb{W}}^{p}( K).
\en

Then, we can obtain the discretized Trefftz-DG  variational formulation corresponding to (\ref{trefftzcontivaria}):
Find $(v_{h},{\bm \sigma}_{h}) \in {\bf V}_h({\cal T}_h)$ such that
\be \label{trefftzvaria}
\mathcal{A}(v_{h},{\bm \sigma}_{h};w,{\bm \tau}) =
\ell(w,{\bm \tau})  \quad \forall (w,{\bm \tau}) \in {\bf V}_h({\cal T}_h).
\en

\section{Error estimates}
In this Section we derive error estimates of approximations generated by the global Trefftz  DG method.

\subsection{The required partition} \label{meshgene}
In order to derive the desired error estimates of the approximate solutions, we require that the partition must satisfy some assumptions. In this part we introduce a kind of particular triangulation
such that these assumptions can be met.

\begin{figure}[H]
\begin{center}
\begin{tabular}{c}
 \epsfxsize=0.5\textwidth\epsffile{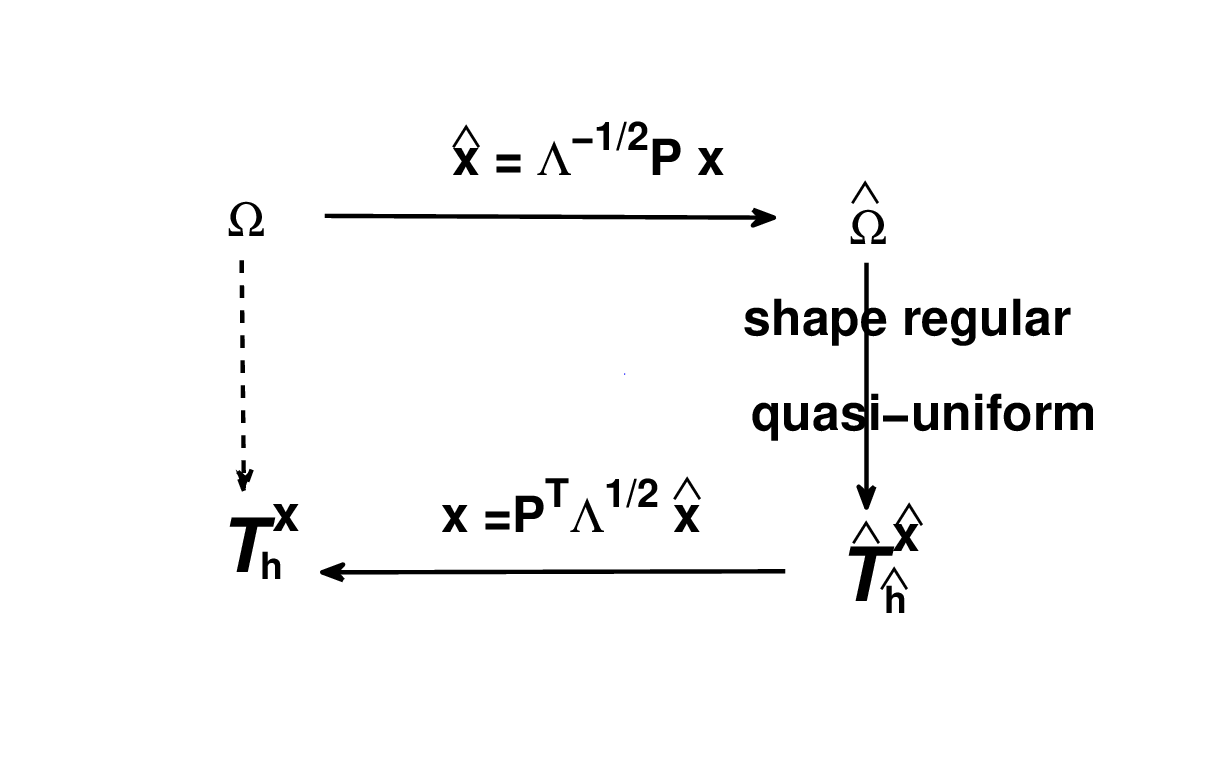}\\
\end{tabular}
\end{center}
 \vskip -0.3in \caption{Mesh generation. } \label{mesh}
\end{figure}

We adopt a non-regularity triangulation ${\cal T}_{h_{{\bf x}}}^{\bf x}$ for the three-dimensional domain \(\Omega\) as follows (see Figure \ref{mesh}).

\noindent {\bf Mesh Generation Algorithm:}

{\it Step 1.}  Determine the transformed domain $\hat\Omega$ under the coordinate transformation (\ref{tran1}).

{\it Step 2.}  Decompose $\hat\Omega$ into polyhedron elements $\{\hat K_{\hat{\bf x}}\}$ such that $\hat{\cal T}_{\hat h_{\hat{\bf x}}}^{\hat{\bf x}}$ is shape regular and quasi-uniform in the usual manner.

{\it Step 3.}  Determine the triangulation ${\cal T}_{h_{{\bf x}}}^{\bf x}=\{K_{\bf x}\}$ of $\Omega$ by using the inverse transformation of (\ref{tran1}) acting on the elements of  $\hat{\cal T}_{\hat h_{\hat{\bf x}}}^{\hat{\bf x}}$.

Under the proposed mesh triangulation, we introduce the two important geometric properties, which are the direct generalization to multidimensional space cases of Lemma 2.2 and 2.3 of \cite{yuan}, respectively.

\begin{lemma} For the proposed triangulation, we have
\be \label{helmrelahhath}
c_0 ||\Lambda^{\frac{1}{2}}||^{-1}h_{{\bf x}} \leq \hat{h}_{\hat{\bf x}} \leq C_0||\Lambda^{\frac{1}{2}}||^{-1}h_{{\bf x}}, ~~\text{and} ~~c_0 ||\Lambda^{\frac{1}{2}}||^{-1}h \leq \hat{h} \leq C_0||\Lambda^{\frac{1}{2}}||^{-1}h,
\en
where $c_0$ and $C_0$ denote two constants independent of $A$.
\end{lemma}

The next Lemma gives a relation between the areas of two bounded hyperplanes based on the coordinate transformation  (\ref{tran1}).

\begin{lemma} \label{arearelanew}
For the proposed triangulation, denote by $\Gamma$ a hyperplane in $\mathbb{R}^d$ which belongs to  $\mathcal{F}_h$, and by $\hat\Gamma$ the correspondingly transformed hyperplane which belongs to $\hat{\mathcal{F}}_{\hat h} $. Then we have
\be \label{arearela}
\frac{|\Gamma|}{|\hat{\Gamma}|}  \leq \text{det}(\Lambda^{\frac{1}{2}})  \lambda_1^{-\frac{1}{2}} ,
\en
where $|f|$ denotes the area of a bounded hyperplane $f$ in the $d$-dimensional space.
\end{lemma}



For the simplicity of notation, let $\rho$ denote the condition number $\text{cond}(A)$ of the anisotropic matrix $A$.
Then $\rho = \text{cond}(\Lambda) = \text{cond}^2(\Lambda^{\frac{1}{2}})= \text{cond}^2(S)$. Without losing generality,  we set $||A||=1$. Equivalently, the original model (\ref{model}) can be normalized such that $||A||=1$.

\subsection{The transformation stability with respect to mesh-dependent norms}

For the global Trefftz DG space ${\bf V}_{h}({\cal T}_{ h})$, 
 we define the following DG norms (see \cite{MP}):
\beq
|||(w,{\bm \tau})|||_{\text{DG}(Q)}^2&= & \frac{1}{2}||c^{-1}\llbracket w\rrbracket_t||^2_{L^2(\mathcal{F}_{h}^{\text{space}})} + \frac{1}{2}||\llbracket {\bm \tau}\rrbracket_t||^2_{{L^2(\mathcal{F}_{h}^{\text{space}})}^d} + \frac{1}{2} ||c^{-1}w||^2_{L^2(\mathcal{F}_{h}^0\cup\mathcal{F}_{h}^T)}
\cr &+& \frac{1}{2}||{\bm \tau}||^2_{L^2(\mathcal{F}_{h}^0\cup\mathcal{F}_{h}^T)^d}
+ ||\alpha^{\frac{1}{2}}\llbracket A^{\frac{1}{4}+\frac{\delta}{2}}w\rrbracket_{\bf N}||^2_{{L^2(\mathcal{F}_{h}^{\text{time}})}^d} + ||\beta^{\frac{1}{2}} \llbracket A^{\frac{1}{2}}{\bm \tau}\rrbracket_{\bf N}||^2_{{L^2(\mathcal{F}_{h}^{\text{time}})}}
\cr & + & \big|\big|\alpha^{\frac{1}{2}} w A^{\frac{1}{4}+\frac{\delta}{2}}{\bf n}^x_{\Omega}\big|\big|^2_{L^2(\mathcal{F}_{h}^D)^d} + ||\beta^{\frac{1}{2}} A^{\frac{1}{2}}{\bm \tau} \cdot {\bf n}^x_{\Omega} ||^2_{{L^2(\mathcal{F}_{h}^N)}},
 \eq
 and
\beq \label{dgplus}
|||(w,{\bm \tau})|||_{\text{DG}(Q)^+}^2 &= & |(w,{\bm \tau})|_{\text{DG}(Q)}^2
  + 2||c^{-1} w^-||^2_{L^2(\mathcal{F}_{h}^{\text{space}})} + 2|| {\bm \tau}^-||^2_{{L^2(\mathcal{F}_{h}^{\text{space}})}^d}
 \cr &+&
 ||\alpha^{-\frac{1}{2}} \{\{ A^{\frac{1}{4}-\frac{\delta}{2}} \bm\tau\}\} ||^2_{{L^2(\mathcal{F}_{h}^{\text{time}})^d}}
 + ||\alpha^{-\frac{1}{2}} A^{\frac{1}{4}-\frac{\delta}{2}} \bm \tau ||^2_{{L^2(\mathcal{F}_{h}^D)^d}}
 \cr &+& ||\beta^{-\frac{1}{2}}\{\{w\}\}||^2_{{L^2(\mathcal{F}_{h}^{\text{time}})}}
+ ||\beta^{-\frac{1}{2}}w||^2_{{L^2(\mathcal{F}_{h}^N)}}. 
 \eq
 In addition, we endow the space $\hat{\bf V}_{\hat h}(\hat{\cal T}_{\hat h})$  with the norm
\beq
|||(\hat w,\hat{\bm \tau})|||_{\text{DG}(\hat Q)}^2&= & \frac{1}{2}||c^{-1}\llbracket \hat w\rrbracket_t||^2_{L^2(\hat{\mathcal{F}}_{\hat h}^{\text{space}})} + \frac{1}{2}||\llbracket \hat {\bm \tau}\rrbracket_t||^2_{{L^2(\hat{\mathcal{F}}_{\hat h}^{\text{space}})}^d} + \frac{1}{2} ||c^{-1}\hat w||^2_{L^2(\hat{\mathcal{F}}_{\hat h}^0\cup\hat{\mathcal{F}}_{\hat h}^T)}
\cr &+& \frac{1}{2}||\hat {\bm \tau}||^2_{L^2(\hat{\mathcal{F}}_{\hat h}^0\cup \hat{\mathcal{F}}_{\hat h}^T)^d}
+ ||\alpha^{\frac{1}{2}} \llbracket \hat w\rrbracket_{\bf N}||^2_{{L^2(\hat{\mathcal{F}}_{\hat h}^{\text{time}})}^d} + ||\beta^{\frac{1}{2}} \llbracket \hat {\bm \tau}\rrbracket_{\bf N}||^2_{{L^2(\hat{\mathcal{F}}_{\hat h}^{\text{time}})}}
\cr & + & ||\alpha^{\frac{1}{2}} \hat w||^2_{L^2(\hat{\mathcal{F}}_{\hat h}^D)} + ||\beta^{\frac{1}{2}} \hat {\bm \tau} \cdot \hat{\bf n}^{\hat x}_{\hat \Omega} ||^2_{{L^2(\hat{\mathcal{F}}_{\hat h}^N)}},
 \eq
 and  the augmented norm
\beq \label{dgplusss}
|||(\hat w,\hat{\bm \tau})|||_{\text{DG}(\hat Q)^+}^2 &= & |(\hat w,\hat{\bm \tau})|_{\text{DG}(\hat Q)}^2
  + 2||c^{-1} \hat w^-||^2_{L^2(\hat{\mathcal{F}}_{\hat h}^{\text{space}})} + 2|| \hat{\bm \tau}^-||^2_{{L^2(\hat{\mathcal{F}}_{\hat h}^{\text{space}})}^d}
 \cr &+&   ||\alpha^{-\frac{1}{2}}
  \{\{\hat{\bm \tau}\}\} 
  ||^2_{{L^2(\hat{\mathcal{F}}_{\hat h}^{\text{time}})^d}}
  +|\alpha^{-\frac{1}{2}}
  \hat{\bm \tau}  ||^2_{{L^2(\hat{\mathcal{F}}_{\hat h}^D)^d}}
  \cr &+& ||\beta^{-\frac{1}{2}}\{\{\hat w\}\}||^2_{{L^2(\hat{\mathcal{F}}_{\hat h}^{\text{time}})}}
| + 
  ||\beta^{-\frac{1}{2}}\hat w||^2_{{L^2(\hat{\mathcal{F}}_{\hat h}^N)}}.
 \eq

The following Lemma states the transformation stability with respect to mesh-dependent norms, which indicates that, in order to obtain the optimal stability estimates and error estimates with respect to $\rho$, the best choice of $\delta$ is set to be $\frac{1}{2}$ throughout the rest paper.

\begin{lemma} \label{helmimportl} For $(w_h,{\bm \tau}_h)\in  {\bf V}_{h}({\cal T}_{ h})$ and $\delta=\frac{1}{2}$, we have
\beq \label{wavestade}
\begin{split}
 |||(w_h,{\bm \tau}_h) |||_{\text{DG}(Q)} & \leq
~\text{det}(\Lambda^{\frac{1}{4}})  \lambda_{\text{min}}^{-\frac{1}{4}}
~ |||(\hat w_{\hat h},\hat{\bm \tau}_{\hat h})|||_{\text{DG}(\hat Q)}, \\
|||(w_h,{\bm \tau}_h)|||_{\text{DG}(Q)^+} & \leq
~\text{det}(\Lambda^{\frac{1}{4}}) \lambda_{\text{min}}^{-\frac{1}{4}}
~ |||(\hat w_{\hat h},\hat{\bm \tau}_{\hat h})|||_{\text{DG}(\hat Q)^+}.
\end{split}
\eq
\end{lemma}

{\it Proof}.  We divide the proof into two steps.

{\it Step 1:} To estimate the terms of $ |||(w_h,{\bm \tau}_h)|||_{\text{DG}(Q)}$ on $\mathcal{F}_{h}^{\text{space}}\cup\mathcal{F}_{h}^0\cup\mathcal{F}_{h}^T$.

By the coordinate transformation (\ref{tran1}) and direct calculation, we obtain
\beq \nonumber
\begin{split}
&
\frac{1}{2}||c^{-1}\llbracket w\rrbracket_t||^2_{L^2(\mathcal{F}_{h}^{\text{space}})} + \frac{1}{2}||\llbracket {\bm \tau}\rrbracket_t||^2_{{L^2(\mathcal{F}_{h}^{\text{space}})}^d} + \frac{1}{2} ||c^{-1}w||^2_{L^2(\mathcal{F}_{h}^0\cup\mathcal{F}_{h}^T)} + \frac{1}{2}||{\bm \tau}||^2_{L^2(\mathcal{F}_{h}^0\cup\mathcal{F}_{h}^T)^d}
 & \\
 & \leq C \text{det}(S^{-1}) ~\bigg(\frac{1}{2}||c^{-1}\llbracket \hat w\rrbracket_t||^2_{L^2(\hat{\mathcal{F}}_{\hat h}^{\text{space}})} + \frac{1}{2}||\llbracket \hat {\bm \tau}\rrbracket_t||^2_{{L^2(\hat{\mathcal{F}}_{\hat h}^{\text{space}})}^d} + \frac{1}{2} ||c^{-1}\hat w||^2_{L^2(\hat{\mathcal{F}}_{\hat h}^0\cup\hat{\mathcal{F}}_{\hat h}^T)}
+ \frac{1}{2}||\hat {\bm \tau}||^2_{L^2(\hat{\mathcal{F}}_{\hat h}^0\cup \hat{\mathcal{F}}_{\hat h}^T)^d}\bigg).
\end{split}
\eq

{\it Step 2:} To estimate the terms of $ |||(w_h,{\bm \tau}_h)|||_{\text{DG}(Q)}$ on $\mathcal{F}_{h}^{\text{time}}\cup\mathcal{F}_{h}^D\cup\mathcal{F}_{h}^N$.
It is easy to check that (here ${\bf n}^t_K=0$)
\be \label{helmne4}
  {\bf n}^{\bf x}_K =
|\Lambda^{\frac{1}{2}}P{\bf n}^{\bf x}_K|~ P^T\Lambda^{\frac{-T}{2}}\hat{\bf n}^{\hat{\bf x}}_{\hat K}.
\en

Combining
\beq \label{vhatvrela}
\begin{split}
&\llbracket A^{\frac{1}{4}+\frac{\delta}{2}}w\rrbracket_{\bf N} = |\Lambda^{\frac{1}{2}}P{\bf n}^{\bf x}_K|~P^T \Lambda^{\frac{\delta}{2}-\frac{1}{4}}\llbracket \hat w\rrbracket_{\bf N},
~~  A^{\frac{1}{4}+\frac{\delta}{2}}{\bf n}^x_{\Omega} = |\Lambda^{\frac{1}{2}}P{\bf n}^{\bf x}_{\Omega}|  P^T \Lambda^{\frac{\delta}{2}-\frac{1}{4}} \hat{\bf n}^{\hat x}_{\hat \Omega},
\\
&
\llbracket A^{\frac{1}{2}}{\bm \tau}\rrbracket_{\bf N} =  |\Lambda^{\frac{1}{2}}P{\bf n}^{\bf x}_K|~\llbracket \hat{\bm \tau}\rrbracket_{\bf N},
~~ A^{\frac{1}{2}}{\bm \tau} \cdot {\bf n}^x_{\Omega} = |\Lambda^{\frac{1}{2}}P{\bf n}^{\bf x}_{\Omega}|  ~\hat {\bm \tau} \cdot \hat{\bf n}^{\hat x}_{\hat \Omega},
\end{split}
\eq
with (\ref{arearela}), 
we get
\beq
\begin{split}
& ||\alpha^{\frac{1}{2}} \llbracket A^{\frac{1}{4}+\frac{\delta}{2}}w\rrbracket_{\bf N}||_{{L^2(\mathcal{F}_{h}^{\text{time}})}^d} \leq \text{det}(\Lambda^{\frac{1}{4}})  \lambda_{\text{min}}^{-\frac{1}{4}} ~ ||\Lambda^{\frac{\delta}{2}-\frac{1}{4}}|| \cdot
||\alpha^{\frac{1}{2}} \llbracket \hat w\rrbracket_{\bf N}||_{{L^2(\hat{\mathcal{F}}_{\hat h}^{\text{time}})}^d},
\\ &
 ||\alpha^{-\frac{1}{2}}
  \{\{A^{\frac{1}{4}-\frac{\delta}{2}} \bm \tau\}\} ||_{{L^2(\mathcal{F}_{h}^{\text{time}})^d}}
  \leq ~ \text{det}(\Lambda^{\frac{1}{4}})  \lambda_{\text{min}}^{-\frac{1}{4}}
  ~ ||\Lambda^{\frac{1}{4}-\frac{\delta}{2}}|| \cdot
 || \{\{\hat{\bm \tau}\}\}||_{{L^2(\hat{\mathcal{F}}_{\hat h}^{\text{time}})^d}},
\\ &
 ||\beta^{\frac{1}{2}} \llbracket A^{\frac{1}{2}}{\bm \tau}\rrbracket_{\bf N}||_{{L^2(\mathcal{F}_{h}^{\text{time}})}} \leq  \text{det}(\Lambda^{\frac{1}{4}})  \lambda_{\text{min}}^{-\frac{1}{4}} ~||\beta^{\frac{1}{2}} \llbracket \hat {\bm \tau}\rrbracket_{\bf N}||_{{L^2(\hat{\mathcal{F}}_{\hat h}^{\text{time}})}},
\\
 &  \big|\big|\alpha^{\frac{1}{2}} w A^{\frac{1}{4}+\frac{\delta}{2}} {\bf n}^x_{\Omega}\big|\big|^2_{L^2(\mathcal{F}_{h}^D)^d}
 \leq \text{det}(\Lambda^{\frac{1}{4}})  \lambda_{\text{min}}^{-\frac{1}{4}}  ~ ||\Lambda^{\frac{\delta}{2}-\frac{1}{4}}|| \cdot
|| \hat w||_{{L^2(\hat{\mathcal{F}}_{\hat h}^D)}},
 \\ & ||\beta^{\frac{1}{2}} A^{\frac{1}{2}}{\bm \tau} \cdot {\bf n}^x_{\Omega} ||^2_{{L^2(\mathcal{F}_{h}^N)}} \leq \text{det}(\Lambda^{\frac{1}{4}})  \lambda_{\text{min}}^{-\frac{1}{4}}  ~
||\beta^{\frac{1}{2}} \hat {\bm \tau} \cdot \hat{\bf n}^{\hat x}_{\hat \Omega} ||_{{L^2(\hat{\mathcal{F}}_{\hat h}^N)}},
  \\ &  ||\beta^{-\frac{1}{2}}\{\{w\}\}||_{{L^2(\mathcal{F}_{h}^{\text{time}})}}
  \leq \text{det}(\Lambda^{\frac{1}{4}})  \lambda_{\text{min}}^{-\frac{1}{4}} ~||\beta^{-\frac{1}{2}}\{\{\hat w\}\}||^2_{{L^2(\hat{\mathcal{F}}_{\hat h}^{\text{time}})}}.
\end{split}
\eq

Combining the two steps with the best choice $\delta=\frac{1}{2}$ yields the desired results (\ref{wavestade}).  \quad $\Box$

\subsection{Error estimates of Trefftz discontinuous Galerkin approximations}
We prove existence and uniqueness of the Trefftz DG solution of (\ref{trefftzvaria}) and the bilinear form in (\ref{adg1}) admits the following upper bounds. Throughout this paper, $C$ denotes a generic positive constant that may have
different values in different occurrences, where $C$ depends on the mesh of $\hat{\cal T}_{\hat h}$ and the shape of the elements.

\begin{lemma} \label{wellposedness}
 There exists a unique solution $(v_{h}, ~{\bm \sigma}_{h})$ to (\ref{trefftzvaria}); moreover, for $\forall (v,{\bm\sigma}), ~(w,{\bm \tau})\in {\bf V}_h({\cal T}_h)$ we have
\be \label{errderi5}
|||(w,{\bm \tau})|||_{\text{DG}(Q)}^2 = \mathcal{A} ((w,{\bm \tau}); (w,{\bm \tau})),
\en
\be \label{errderi6}
|\mathcal{A} (v,{\bm\sigma}; w,{\bm \tau})|\leq
|||(v,{\bm\sigma})|||_{\text{DG}(Q)^+}~|||(w,{\bm \tau})|||_{\text{DG}(Q)},
\en
and
\be \label{errderi7}
|\mathcal{A} (v,{\bm\sigma}; w,{\bm \tau})|\leq 2 |||(v,{\bm\sigma})|||_{\text{DG}(Q)} ~|||(w,{\bm \tau})|||_{\text{DG}(Q)^+}.
\en
 \end{lemma}

{\it Proof}.  Provided that $(v,{\bm\sigma}), ~(w,{\bm \tau}) \in {\bf V}_h({\cal T}_h)$, local integration by parts permits us to rewrite the
bilinear form $\mathcal{A} ((v,{\bm\sigma}), ~(w,{\bm \tau}))$ as

\begin{eqnarray} \label{adg2}
&\mathcal{A}(v,{\bm\sigma}; w,{\bm \tau})=  -\int_{\mathcal{F}_h^{\text{space}}}\big(c^{-2}\llbracket v\rrbracket_t w^+ + \llbracket {\bm \sigma} \rrbracket_t \cdot {\bm \tau}^+ \big)~d{\bf x}
 +  \int_{\mathcal{F}_h^0} \big( c^{-2}v w + {\bm \sigma} \cdot {\bm \tau}\big)~d{\bf x}
\cr &+ \int_{\mathcal{F}_h^{\text{time}}}
\big( -\llbracket A^{\frac{1}{2}}v \rrbracket_{\bf N} \cdot \{\{{\bm \tau}\}\}
-\llbracket A^{\frac{1}{2}}{\bm \sigma} \rrbracket_{\bf N} \{\{ w\} \}
 + \alpha\llbracket A^{\delta}v \rrbracket_{\bf N} \cdot\llbracket A^{\frac{1}{2}}w \rrbracket_{\bf N}
  + \beta \llbracket A^{\frac{1}{2}}{\bm \sigma} \rrbracket_{\bf N} \llbracket A^{\frac{1}{2}}{\bm \tau} \rrbracket_{\bf N}  \big)~dS
      \cr & +   \int_{\mathcal{F}_{h}^D} ( -v A^{\frac{1}{2}}{\bf n}^x_{\Omega}\cdot{\bm \tau}    + \alpha v w A^{\delta}{\bf n}^x_{\Omega} \cdot A^{\frac{1}{2}}{\bf n}^x_{\Omega}) ~dS
      \cr & + \int_{\mathcal{F}_{h}^N} (  -(A^{\frac{1}{2}}{\bm \sigma}\cdot {\bf n}^x_{\Omega})w + \beta  (A^{\frac{1}{2}}{\bm \sigma} \cdot {\bf n}^x_{\Omega})
  (A^{\frac{1}{2}}{\bm \tau}\cdot {\bf n}^x_{\Omega}) )~dS.
 \end{eqnarray}

By taking $(v,{\bm\sigma})=(w,{\bm \tau})$ and summing the two expressions given in (\ref{adg1}) and (\ref{adg2}), we obtain (\ref{errderi5}).

$|||\cdot|||_{\text{DG}(Q)}$ and $|||\cdot|||_{\text{DG}(Q)^+}$ are only seminorms on broken Sobolev spaces defined on the mesh ${\cal T}_h(Q)$, but are
norms on ${\bf T}({\cal T}_h(Q))$: indeed $|||(w,{\bm\tau})|||_{\text{DG}}=0$
 for $(w,{\bm\tau})\in {\bf T}({\cal T}_h(Q))$ implies that $(w,{\bm\tau})$ is
solution of the homogeneous IBVP (\ref{model}) with zero initial and boundary conditions, so $(w,{\bm\tau})={\bf 0}$ by the well-posedness of the IBVP itself (see \cite[Section 5.1]{MP} and \cite[Lemma 4.1]{Kretzschmar}). Thus the variational formulation (\ref{trefftzcontivaria}) has a unique discrete solution. Existence of the solution follows from linearity of the problem and finite dimensionality.

By applying the Cauchy-Schwarz inequality to (\ref{adg1}) and (\ref{adg2}), respectively,
we obatin (\ref{errderi6}) and (\ref{errderi7}). \quad $\Box$

\begin{theorem} \label{anisoerr}
Assume that the IBVP solution $(v,{\bm \sigma})\in C^{k_t-1}(I;H^{k_{\bf x}+1}(\Omega))\times C^{k_t}(I;H^{k_{\bf x}}(\Omega)^d)$, and that $s=\text{min}\{p,k_t-1,k_{\bf x}-1\}$. Then we have, 
\beq\nonumber
\frac{1}{2}\bigg(||c^{-1}  (v-v_h)||_{L^2(\Omega\times\{T\})} + ||\bm\sigma-\bm\sigma_h||_{L^2(\Omega\times\{T\})^d}\bigg)
 \leq  |||(v,{\bm \sigma}) - (v_h,{\bm \sigma}_{h})|||_{\text{DG}(Q)}
 \leq C\rho^{\frac{1}{4}}h^{s+\frac{1}{2}}
 |(v,{\bm \sigma})|_{H^{s+1}(Q)^{1+d}}.
\eq
\end{theorem}

{\it Proof}. By (\ref{errderi5}), (\ref{trefftzcontivaria}) and (\ref{trefftzvaria}), we obtain, for $\forall (w_h,{\bm \tau}_h)\in {\bf V}_h({\cal T}_h(Q))$,

\be \label{inteerr20}
|||(v,{\bm \sigma}) - (v_h,{\bm \sigma}_h)|||_{\text{DG}(Q)}^2 = \mathcal{A}_{\text{DG}(Q)} ((v,{\bm \sigma}) - (v_h,{\bm \sigma}_h); (v,{\bm \sigma})-(w_h,{\bm \tau}_h)).
\en
Taking into account (\ref{errderi7}), we get the abstract error estimate:
\be \label{abstarct}
|||(v,{\bm \sigma}) - (v_h,{\bm \sigma}_h)|||_{\text{DG}(Q)} \leq
2\mathop{\text{inf}}\limits_{ (w_h,{\bm \tau}_h) \in {\bf V}_h({\cal T}_h(Q))
} |||(v,{\bm \sigma}) - (w_h,{\bm \tau}_h)|||_{\text{DG}(Q)^+}.
\en

By the existing approximation result \cite[Corollary 4]{MP}, there exists $\hat Q_{\hat h}(\hat v,\hat{\bm\sigma})=(\hat Q_{\hat h}\hat v,\hat Q_{\hat h}\hat{\bm\sigma})\in \hat{\bf V}_{\hat h}(\hat{\cal T}_{\hat h}(\hat Q))$ such that,
\be \label{isotreeappr4}
|(\hat v,\hat{\bm\sigma})-\hat Q_{\hat h}(\hat v,\hat{\bm\sigma})|_{H^j(\hat K)^{1+d}} \leq C \hat h^{s+1-j}|(\hat v,\hat{\bm\sigma})|_{H^{s+1}(\hat K)^{1+d}}.
\en
Using the inverse transformation of (\ref{firstorderrela}), set
\be \label{transvari}
Q_h(v,\bm\sigma)= (\hat Q_{\hat h}\hat v, P^T\hat Q_{\hat h}\hat{\bm\sigma}).
\en

Using (\ref{abstarct}), (\ref{transvari}) and (\ref{wavestade}), we obtain
\begin{eqnarray} \label{inteapp21}
& |||(v,{\bm \sigma}) - (v_h,{\bm \sigma}_h)|||_{\text{DG}(Q)} \leq
|||(v,{\bm \sigma}) - Q_h(v,\bm\sigma)|||_{\text{DG}(Q)^+}
\cr & \leq
~2\text{det}(\Lambda^{\frac{1}{4}})  \lambda_{\text{min}}^{-\frac{1}{4}}
~ |||(\hat v,\hat{\bm\sigma})-\hat Q_{\hat h}(\hat v,\hat{\bm\sigma})|||_{\text{DG}(\hat Q)^+}.
\end{eqnarray}

By the transformation stability of Lemma \ref{helmimportl}, we only need to derive a bound of $(\hat v,\hat{\bm\sigma})-\hat Q_{\hat h}(\hat v,\hat{\bm\sigma})$ in terms of elementwise sums of traces, tracking the dependence on spatial and temporal meshsizes. Taking into account the defintion (\ref{dgplusss}), we get, for $(\hat w,\hat{\bm \tau}) \in \hat{\bf T}(\hat{\cal T}_{\hat h}(\hat Q))$,
\beq \label{errderi8}
|||(\hat w,\hat{\bm \tau})|||_{\text{DG}(\hat Q)^+}^2 & \leq & \sum\limits_{\hat K=\hat K_{\hat x}\times I_{n'}\in \hat{\cal T}_{\hat h}(\hat Q)}\bigg[ ||c^{-1}\hat w||^2_{L^2(\hat K_{\hat x}\times \{t_{n'-1},t_{n'}\})} + ||\hat{\bm \tau}||^2_{L^2(\hat K_{\hat x}\times \{t_{n'-1},t_{n'}\})^d}
\cr & + & \sum\limits_{F\in \partial \hat K\cap \hat Q\cap (\hat{\mathcal{F}}_{\hat h}^{\text{time}}\cup \hat{\mathcal{F}}_{\hat h}^D)} ||\hat w \hat{\bf n}^{\hat x}_{\hat F}||_{L^2(\hat F)^d}^2 + \sum\limits_{\hat F\in \partial \hat K\cap \hat Q\cap (\hat{\mathcal{F}}_{\hat h}^{\text{time}}\cup \hat{\mathcal{F}}_{\hat h}^N)} ||\hat w||_{L^2(\hat F)}^2
\cr & + & \sum\limits_{\hat F\in \partial \hat K\cap \hat Q\cap (\hat{\mathcal{F}}_{\hat h}^{\text{time}}\cup \hat{\mathcal{F}}_{\hat h}^N)} ||\hat {\bm \tau}\cdot \hat{\bf n}^{\hat x}_{\hat F}||_{L^2(\hat F)}^2 + \sum\limits_{\hat F\in \partial \hat K\cap \hat Q\cap (\hat{\mathcal{F}}_{\hat h}^{\text{time}}\cup \hat{\mathcal{F}}_{\hat h}^D)} ||\hat{\bm \tau}||_{L^2(\hat F)^d}^2
\bigg].
 \eq
By the standard weighted trace inequality applied in the time and space directions independently (see \cite[Sec 1.6.6]{brenner}), the following bound holds true:
\beq \label{errderi9}
|||(\hat w,\hat{\bm \tau})|||_{\text{DG}(\hat Q)^+}^2 & \leq & C \sum\limits_{\hat K=\hat K_{\hat x}\times I_{n'}\in \hat{\cal T}_{\hat h}(\hat Q)}\bigg[ \hat h^{-1}_{n'}\bigg(||c^{-1}\hat w||^2_{L^2(\hat K)}
+ ||\hat{\bm \tau}||^2_{L^2(\hat K)^d} \bigg)
\cr &+&  \hat h_{n'}\bigg(|c^{-1}\hat w|^2_{H^1(I_{n'};L^2(\hat K_{\hat x}))} + |\hat{\bm \tau}|^2_{H^1(I_{n'};L^2(\hat K_{\hat x})^d)} \bigg)
\cr & + & \hat h_{\hat K_{\hat x}}^{-1} \bigg( ||\hat w||^2_{L^2(\hat K)} + ||\hat w||^2_{L^2(\hat K)} \bigg)
+ \hat h_{\hat K_{\hat x}}  \bigg( | \hat w|^2_{L^2 (I_{n'};H^1(\hat K_{\hat x}))} + |\hat w|^2_{L^2 (I_{n'};H^1(\hat K_{\hat x}))}  \bigg)
\cr & + &  \hat h_{\hat K_{\hat x}}^{-1} ||\hat{\bm \tau}||^2_{L^2(\hat K)^d} +
\hat h_{\hat K_{\hat x}} |\hat{\bm \tau}|^2_{L^2 (I_{n'};H^1(\hat K_{\hat x})^d)}
\bigg].
 \eq

Using (\ref{inteapp21}), (\ref{errderi9}), (\ref{isotreeappr4}), (\ref{helmrelahhath}) and the scaling argument, we obtain
 \beq \label{errderi10} 
 &|||(v,{\bm \sigma}) - (v_h,{\bm \sigma}_h)|||_{\text{DG}(Q)}
  \leq ~C\text{det}(\Lambda^{\frac{1}{4}})  \lambda_{\text{min}}^{-\frac{1}{4}}
    \sum\limits_{\hat K=\hat K_{\hat x}\times I_{n'}\in \hat {\cal T}_{\hat h}(\hat Q)}\bigg[
\hat h^{-\frac{1}{2}}||(\hat v,\hat {\bm \sigma})-\hat Q_{\hat h}(\hat v,\hat {\bm \sigma})||_{L^2(\hat  K)^{1+d}}
\cr  &+  \hat h^{\frac{1}{2}}\big|\big((\hat v,\hat {\bm \sigma})-\hat Q_{\hat h}(\hat v,\hat {\bm \sigma})\big)\big|_{H^1(\hat K)^{1+d}}\bigg]
 \leq  C \text{det}(\Lambda^{\frac{1}{4}}) ~\lambda_{\text{min}}^{-\frac{1}{4}}  ~\hat h^{s+\frac{1}{2}} \sum\limits_{\hat K=\hat K_{\hat x}\times I_{n'}\in \hat {\cal T}_{\hat h}(\hat Q)}|(\hat v,\hat {\bm \sigma})|_{H^{s+1}(\hat K)^{1+d}}
\cr 
 &\leq C\rho^{\frac{1}{4}}~h^{s+{\frac{1}{2}}}
|(v,{\bm \sigma})|_{H^{s+1}(Q)^{1+d}}. \quad\Box
 \eq

\section{A standard Trefftz DG method}
A natural idea is to apply the standard Trefftz DG method to the isotropic wave equation (\ref{transmodel}) derived by the coordinate transformation $S$ and the scaled transformation (\ref{firstorderrela}), and then use the image of the resulting approximation under the inverse transformation $S^{-1}$ and the inverse scaled transformation (\ref{firstorderrela}) as the desired approximation of $(v, {\bm\sigma})$. We will give the detailed derivation of variational formulation in this Section.

Let $(\hat v, \hat{\bm\sigma})$ 
denote the analytic solution of the equation (\ref{transmodel}) with the transformed boundary  and initial conditions from the original boundary and initial conditions (\ref{model}):
  \begin{equation} \label{transhomomodel}
\left\{ \begin{aligned}
      &  \hat v(\cdot,0)=v_0, \quad \hat{\bm \sigma}(\cdot,0)= P{\bm \sigma}_0 & \text{on} \quad \hat\Omega, \\
       &    \hat v = g_D & \text{on} \quad \hat\Gamma_D\times[0,T],
       \\
      &   |\Lambda^{\frac{1}{2}}P{\bf n}^{\bf x}_{\Omega}| \hat{\bm \sigma}\cdot \hat{\bf n}^x_{\hat\Omega} = g_N & \text{on} \quad \hat\Gamma_N\times[0,T].
                          \end{aligned} \right.
                          \end{equation}

Define the numerical fluxes as follows.
\begin{equation} \label{isostandflux} \nonumber
\breve{v}=\left\{ \begin{aligned}
     &  \hat v^- \\
     & \hat v \\
     & v_0 \\
     & \{\{\hat v\}\} + \beta \llbracket  \hat{\bm\sigma} \rrbracket_N \\
     & g_D \\
     & \hat v + \beta ( |\Lambda^{\frac{1}{2}}P{\bf n}^{\bf x}_{\Omega}| \hat{\bm \sigma}\cdot \hat{\bf n}^x_{\hat\Omega} - g_N)
  \end{aligned} \right.
  \breve{\bm \sigma}= \left\{ \begin{aligned}
     &  \hat{\bm \sigma}^- &\quad {\text on} ~ \hat{\mathcal{F}}_{\hat h}^{\text{space}},\\
     & \hat{\bm \sigma} &\quad {\text on} ~ \hat{\mathcal{F}}_{\hat h}^T, \\
     & P{\bm \sigma}_0 &\quad {\text on} ~ \hat{\mathcal{F}}_{\hat h}^0, \\
     &  \{\{\hat{\bm \sigma}\}\} + \alpha \llbracket  \hat v \rrbracket_N &\quad {\text on} ~ \hat{\mathcal{F}}_{\hat h}^{\text{time}}, \\
     & \hat{\bm \sigma} + \alpha (\hat v-g_D)\hat{\bf n}^x_{\hat\Omega} &\quad {\text on} ~ \hat{\mathcal{F}}_{\hat h}^D, \\
     & g_N |\Lambda^{\frac{1}{2}}P{\bf n}^{\bf x}_{\Omega}|^{-1}\hat{\bf n}^x_{\hat\Omega}   &\quad {\text on} ~ \hat{\mathcal{F}}_{\hat h}^N.
  \end{aligned} \right.
                          \end{equation}
Then \((\hat v, \hat{\bm\sigma})\in  \hat{\bf T}(\hat{\cal T}_{\hat h})\) satisfies \be \label{isotrefftzcontivaria}
\hat{\mathcal{A}}(\hat v,\hat{\bm \sigma};\hat w,\hat{\bm \tau}) =
\hat{\ell}(\hat w, \hat{\bm \tau})  \quad \forall (\hat w,\hat{\bm \tau}) \in \hat{\bf T}(\hat{\cal T}_{\hat h}),
\en
where
\beq \label{isoadg1}
\hat{\mathcal{A}}(v,{\bm \sigma};w,{\bm \tau})&= & \int_{\hat{\mathcal{F}}_{\hat h}^{\text{space}}}\big(c^{-2}v^-\llbracket w\rrbracket_t + {\bm \sigma}^-\cdot \llbracket {\bm \tau} \rrbracket_t  \big)~d\hat{\bf x}
+  \int_{\hat{\mathcal{F}}_{\hat h}^{\text{time}}}
\bigg(\{\{v\}\} \llbracket {\bm \tau} \rrbracket_{\bf N} +
 \{\{ {\bm \sigma} \}\} \cdot \llbracket w \rrbracket_{\bf N}
 \cr &+&  \alpha \llbracket v \rrbracket_{\bf N} \cdot \llbracket w \rrbracket_{\bf N}
  + \beta \llbracket {\bm \sigma} \rrbracket_{\bf N} \llbracket {\bm \tau} \rrbracket_{\bf N}  \bigg)~dS
    \cr & + &  \int_{\hat{\mathcal{F}}_{\hat h}^T} \big( c^{-2}v w + {\bm \sigma} \cdot {\bm \tau}\big)~d\hat{\bf x}
    +  \int_{\hat{\mathcal{F}}_{\hat h}^D} ({\bm \sigma} \cdot \hat{\bf n}^x_{\hat\Omega} w + \alpha v w ) ~dS
    \cr & + &   \int_{\hat{\mathcal{F}}_{\hat h}^N} \bigg(  v ({\bm \tau}\cdot \hat{\bf n}^x_{\hat\Omega} ) + \beta  |\Lambda^{\frac{1}{2}}P{\bf n}^{\bf x}_{\Omega}| ({\bm \sigma} \cdot \hat{\bf n}^x_{\hat\Omega} )
  ({\bm \tau}\cdot \hat{\bf n}^x_{\hat\Omega} ) \bigg)~dS,
 \eq
 and
\beq
 \hat{\ell}(w,{\bm \tau}) &=& \int_{\hat{\mathcal{F}}_{\hat h}^0}\big( c^{-2}v_0 w + P{\bm \sigma}_0 \cdot {\bm \tau}\big)~d\hat{\bf x}
+   \int_{\hat{\mathcal{F}}_{\hat h}^N} g_N \bigg( \beta {\bm \tau}\cdot \hat{\bf n}^x_{\hat\Omega}  - |\Lambda^{\frac{1}{2}}P{\bf n}^{\bf x}_{\Omega}|^{-1} w \bigg) ~ dS
    \cr & + &  \int_{\hat{\mathcal{F}}_{\hat h}^D} \alpha g_D w ~ dS
 - \int_{\hat{\mathcal{F}}_{\hat h}^D} g_D {\bm \tau}\cdot \hat{\bf n}^x_{\hat\Omega} ~ dS.
 \eq

\begin{remark}
The choice of numerical fluxes $(\breve{v}, \breve{\bm \sigma})$
defined on ``time-like" interfaces and boundary faces is such that
the Trefftz DG formulation is consistent; namely, if $(\hat v,
\hat{\bm\sigma})\in H^1(\hat Q)$ solves (\ref{transmodel}) with boundary and initial conditions (\ref{transhomomodel}), then it satisfies (\ref{isotrefftzcontivaria}). In particular, the flux $\breve{\bm \sigma}$ satisfies $|\Lambda^{\frac{1}{2}}P{\bf n}^{\bf x}_{\Omega}| \breve{\bm \sigma}\cdot \hat{\bf n}^x_{\hat\Omega} = g_N$
 coinciding with the Neumann boundary condition satisfied by the exact solution $(\hat v,\hat{\bm\sigma})$.
\end{remark}

Denote by $(\hat v_{\hat h}, \hat{\bm\sigma}_{\hat h})$ the discrete approximation of  $(\hat v, \hat{\bm\sigma})$, and let $(\tilde{v}_h({\bf x}), \tilde{\bm\sigma}_h({\bf x}))$ denote the image of $(\hat v_{\hat h}, \hat{\bm\sigma}_{\hat h})$ under the inverse scaled transformation of (\ref{firstorderrela}) and the coordinate transformation (\ref{tran1}). By the definition (\ref{adg1}) of the sesquilinear form $\mathcal{A}(\cdot; \cdot)$, the definition (\ref{isoadg1}) of $\hat{\mathcal{A}}(\cdot; \cdot)$, (\ref{helmne4}), and (\ref{vhatvrela}), we have
\begin{equation} \label{helmpwdgrelafhb}
\left\{ \begin{aligned}
& c^{-2}v^-\llbracket w\rrbracket_t &=& \quad
c^{-2} \hat v^-\llbracket \hat w\rrbracket_t \\
&  {\bm \sigma}^-\cdot \llbracket {\bm \tau} \rrbracket_t &=& \quad
\hat{\bm \sigma}^-\cdot \llbracket \hat{\bm \tau} \rrbracket_t
   \end{aligned} \right.   \quad {\text on} ~ \mathcal{F}_h^{\text{space}},
                          \end{equation}
and
\begin{equation} \label{helmpwdgrelafhi}
\left\{ \begin{aligned}
& \{\{v\}\} \llbracket {A^{\frac{1}{2}}\bm \tau} \rrbracket_{\bf N}
&=& \quad  |\Lambda^{\frac{1}{2}}P{\bf n}^{\bf x}_K|~ \{\{\hat v\}\}   \llbracket \hat{\bm \tau}\rrbracket_{\bf N} \\
 &\{\{ {\bm \sigma} \}\} \cdot \llbracket A^{\frac{1}{2}}w \rrbracket_{\bf N}
 & = &\quad |\Lambda^{\frac{1}{2}}P{\bf n}^{\bf x}_K|~\{\{\hat{\bm \sigma}\}\} ~\llbracket \hat w\rrbracket_{\bf N}\\
  & \alpha \llbracket A^{\frac{1}{2}}v \rrbracket_{\bf N} \cdot \llbracket A^{\frac{1}{2}}w \rrbracket_{\bf N} & = &
  \quad \alpha |\Lambda^{\frac{1}{2}}P{\bf n}^{\bf x}_K|^2~ \llbracket \hat v \rrbracket_{\bf N} \cdot \llbracket \hat w \rrbracket_{\bf N} \\
   & \beta \llbracket A^{\frac{1}{2}}{\bm \sigma} \rrbracket_{\bf N} \llbracket A^{\frac{1}{2}}{\bm \tau} \rrbracket_{\bf N}   &=&
   \quad \beta   |\Lambda^{\frac{1}{2}}P{\bf n}^{\bf x}_K|^2~ \llbracket \hat{\bm \sigma} \rrbracket_{\bf N} \llbracket \hat{\bm \tau} \rrbracket_{\bf N}
                          \end{aligned} \right.   \quad {\text on}  ~ \mathcal{F}_h^{\text{time}}.
                          \end{equation}
Thus we can see that, no matter how the stabilization parameters $\alpha$ and $\beta$ defined on $\mathcal{F}_h^{\text{time}}$ are chosen, the sesquilinear form $\mathcal{A}(\cdot; \cdot)$ for the anisotropic case can not coincide with $\mathcal{\hat A}_{\hat h}(\cdot; \cdot)$ in the sense of proportionality. Thus the proposed approximation $({v}_h({\bf x}), {\bm\sigma}_h({\bf x}))$ is different from the image $(\tilde{v}_h({\bf x}), \tilde{\bm\sigma}_h({\bf x}))$. However, we still build the desired error estimates of $(\tilde{v}_h({\bf x}), \tilde{\bm\sigma}_h({\bf x}))$ as follows.

\begin{theorem} \label{isoerr}
Assume that the IBVP solution $(v,{\bm \sigma})\in C^{k_t-1}(I;H^{k_{\bf x}+1}(\Omega))\times C^{k_t}(I;H^{k_{\bf x}}(\Omega)^d)$, and that $s=\text{min}\{p,k_t-1,k_{\bf x}-1\}$. Then we have, 
\beq\nonumber
& \frac{1}{2}&||c^{-1}  (v-\tilde{v}_h)||_{L^2(\Omega\times T)} + \frac{1}{2}||\bm\sigma-\tilde{\bm\sigma}_h||_{L^2(\Omega\times T)^d}
 \leq  |||(v,{\bm \sigma}) - (\tilde v_h,\tilde{\bm \sigma}_{h})|||_{\text{DG}(Q)}
 \cr &\leq&  C\rho^{\frac{1}{4}}~h^{s+\frac{1}{2}}
 |(v,{\bm \sigma})|_{H^{s+1}(Q)^{1+d}}.
\eq
\end{theorem}

{\it Proof}. By Lemma \ref{helmimportl} and the scaling argument, we have
\be \label{isofirstapp10}
|||(v,{\bm \sigma}) - ( \tilde v_{h},\tilde{\bm \sigma}_{h})|||_{\text{DG}(Q)}
\leq  \text{det}(\Lambda^{\frac{1}{4}}) ~\lambda_{\text{min}}^{-\frac{1}{4}}
 |||(\hat v,\hat{\bm \sigma}) - (\hat v_{\hat h},\hat{\bm \sigma}_{\hat h})|||_{\text{DG}(\hat Q)}.
\en

By the abstract error estimate (\ref{abstarct}), 
approximation result (\ref{isotreeappr4}) and trace inequalities (\ref{errderi8})-(\ref{errderi10}) for the isotropic case, we get
\begin{eqnarray}\label{isofirstapp100}
& |||(v,{\bm \sigma}) - ( \tilde v_{h},\tilde{\bm \sigma}_{h})|||_{\text{DG}(Q)}
 \leq  \text{det}(\Lambda^{\frac{1}{4}}) ~\lambda_{\text{min}}^{-\frac{1}{4}} ~\hat h^{s+\frac{1}{2}} \sum\limits_{\hat K=\hat K_{\hat x}\times I_{n'}\in \hat {\cal T}_{\hat h}(\hat Q)}|(\hat v,\hat {\bm \sigma})|_{H^{s+1}(\hat K)^{1+d}}
\cr & \leq C\rho^{\frac{1}{4}}~h^{s+{\frac{1}{2}}}
|(v,{\bm \sigma})|_{H^{s+1}(Q)^{1+d}}. \quad \Box
 \end{eqnarray}

\begin{remark}
By Theorems \ref{anisoerr} and \ref{isoerr}, we can see that the proposed method and the standard Trefftz DG method have the same convergence order with respect to $h$ and $\rho$, and  almost have the same computational cost. Besides, we believe that the orders of the condition number $\rho$ in the error estimates are optimal
since the transformation stability estimates seem sharp.
\end{remark}

\begin{remark}\label{remarkcomboundhomo}
We address that, compared with the valid error analysis of \cite{BMPS} only for the homogeneous Neumann boundary conditions from the original IBVP and part valid numerical results for nonhomogeneous Neumann boundary conditions, there is no constraint for our proposed methods introduced in Sections 5 and 6 on the data of homogeneous Neumann boundary conditions from the original IBVP, owing to the use of trace estimate in (\ref{errderi9}) instead of  the employed inverse trace estimate as in \cite[Section 6.1]{BMPS}.
\end{remark}

\section{A nonhomogeneous model}
The model reads as
  \begin{equation} \label{nonmodel}
\left\{ \begin{aligned}
     &   A^{\frac{1}{2}}\nabla v + \frac{\partial {\bm \sigma}}{\partial t}={\bf 0} & \text{in} \quad Q, \\
      & \nabla\cdot A^{\frac{1}{2}}{\bm \sigma} + c^{-2} \frac{\partial v}{\partial t}=f & \text{in} \quad Q, \\
      &  v(\cdot,0)=v_0, \quad {\bm \sigma}(\cdot,0)= {\bm \sigma}_0 & \text{on} \quad \Omega, \\
       &    v = g_D & \text{on} \quad \Gamma_D\times[0,T],
       \\
      &   A^{\frac{1}{2}}{\bm \sigma} \cdot {\bf n}^x_{\Omega} = g_N & \text{on} \quad \Gamma_N\times[0,T].
                          \end{aligned} \right.
                          \end{equation}

In the framework of the global Trefftz DG method combined with overlapping local DG method, we decompose the solution $(v, ~{\bm \sigma})$
of the problem \((\ref{nonmodel})\) into $(v, ~ {\bm \sigma})= (v^{(1)}, ~{\bm \sigma}^{(1)}) + (v^{(2)}, ~{\bm \sigma}^{(2)})$,
where \((v^{(1)}, ~{\bm \sigma}^{(1)})\) is a particular local solution of the first two equations (\ref{nonmodel}) on each fictitious domain with homogeneous boundary and initial conditions, and \((v^{(2)}, ~{\bm \sigma}^{(2)})\) satisfies the locally homogeneous wave equation.

Similarly to the derivation of (\ref{trefftzcontivaria}), we can obtain the Trefftz-DG  variational formulation:
Find $(v^{(2)},{\bm \sigma}^{(2)}) \in {\bf T}({\cal T}_h)$ such that
\be \label{nontrefftzcontivaria}
\mathcal{A}(v^{(2)},{\bm \sigma}^{(2)};w,{\bm \tau}) =
\tilde{\ell}(w,{\bm \tau})- \mathcal{A}(v^{(1)},{\bm \sigma}^{(1)};w,{\bm \tau})  \quad \forall (w,{\bm \tau}) \in {\bf T}({\cal T}_h),
\en
where $\mathcal{A}(\cdot;\cdot)$ is defined as (\ref{adg1}), and
\beq \label{adgtildelf1}
 \tilde{\ell}(w,{\bm \tau}) &=& \int_Q f w ~dV
 + \int_{\mathcal{F}_h^0}\big( c^{-2}v_0 w + {\bm \sigma}_0 \cdot {\bm \tau}\big)~d{\bf x}
    \cr & + &  \int_{\mathcal{F}_h^D} \alpha g_D w A^{\delta}{\bf n}^x_{\Omega}\cdot A^{\frac{1}{2}}{\bf n}^x_{\Omega} ~ dS
 - \int_{\mathcal{F}_h^D} g_D A^{\frac{1}{2}}{\bm \tau}\cdot {\bf n}^x_{\Omega} ~ dS
   \cr & + &  \int_{\mathcal{F}_h^N} g_N \bigg( \beta A^{\frac{1}{2}}{\bm \tau}\cdot {\bf n}^x_{\Omega}  - w \bigg) ~ dS.
 \eq

 \subsection{Nonhomogeneous local problems} For each space-time element $K=K_{\bf x}\times I_n\in {\cal T}_h, K_{\bf x}\in {\cal T}_{h_{{\bf x}}}^{\bf x}$, let $K_{\bf x}^{\ast}$ be a fictitious domain that 
contains $K_{\bf x}$ as its subdomain. Set the fictitious domain $K^{\ast}=K_{\bf x}^{\ast}\times I_n$, $\mathcal{F}_{K^{\ast}}^{t_n}=K_{\bf x}^{\ast}\times \{t=t_n\}$, and $\mathcal{F}_{K^{\ast},n}^D=\partial K_{\bf x}^{\ast}\times I_n$.

 Define ${\bf V}^{(1)}_{K^{\ast}}=H^1(I_n;H^{2,2}(K_{\bf x}^{\ast}))\times H^1(I_n;H^{1,1}(K_{\bf x}^{\ast})^d)$. The particular solution
$(v^{(1)}, ~{\bm \sigma}^{(1)}) \in (L^2(Q))^{1+d}$ is defined as
$(v^{(1)},~ {\bm \sigma}^{(1)}) \big|_{K} =(v^{(1)}_K, ~{\bm \sigma}^{(1)}_K)=(v^{(1)}_{K^{\ast}}, ~{\bm \sigma}^{(1)}_{K^{\ast}})\big|_{K}$,
where $(v^{(1)}_{K^{\ast}}, ~{\bm \sigma}^{(1)}_{K^{\ast}}) \in {\bf V}^{(1)}_{K^{\ast}}$ 
 satisfies the nonhomogeneous {\it local} acoustic equation on the fictitious domain $K^{\ast}$:
\begin{equation} \label{localnonhomomacro}
\left\{ \begin{aligned}
    &   A^{\frac{1}{2}}\nabla v_{K^{\ast}}^{(1)} + \frac{\partial {\bm \sigma}_{K^{\ast}}^{(1)}}{\partial t}={\bf 0} & \text{in} \quad K^{\ast}, \\
      & \nabla \cdot A^{\frac{1}{2}}{\bm \sigma}_{K^{\ast}}^{(1)} + c^{-2} \frac{\partial v_{K^{\ast}}^{(1)}}{\partial t}=f & \text{in} \quad K^{\ast}, \\
      &  v_{K^{\ast}}^{(1)} =0,~~ {\bm \sigma}_{K^{\ast}}^{(1)}={\bf 0} & \quad\quad\quad \text{on} \quad K_{\bf x}^{\ast} \times\{t=t_{n-1}\}, \hfill\\
        &    v_{K^{\ast}}^{(1)} = 0  & \text{on} \quad \partial K_{\bf x}^{\ast}\times I_n .
       \\
                    \end{aligned} \right.
                          \end{equation}

Similarly to the derivation of (\ref{stavaria0})-(\ref{stavaria}), we can obtain
\beq \label{localvari}
 -\int_{K^{\ast}} && \bigg( v(\nabla\cdot A^{\frac{1}{2}} {\bm \tau} + c^{-2}\frac{\partial w}{\partial t}) + {\bm \sigma} \cdot(A^{\frac{1}{2}}\nabla w + \frac{\partial \bm\tau}{\partial t})\bigg)dV
 +\int_{\partial {K^{\ast}}}\bigg( \check{v}( A^{\frac{1}{2}} {\bm \tau}\cdot {\bf n}^x_K + c^{-2} w n^t_K)
 \cr +& & \check{\bm \sigma} \cdot (w A^{\frac{1}{2}}{\bf n}^x_K + {\bm \tau} n^t_K) \bigg)dS = \int_{K^{\ast}} f~w~dV, ~\forall (w,{\bm \tau}) \in {\bf V}^{(1)}_{K^{\ast}}.
\eq
Define the numerical fluxes as follows.

\begin{equation} \label{standflux1} \nonumber
\check{v}=\left\{ \begin{aligned}
     & v \\
     & 0 \\
     & 0 \\
  \end{aligned} \right. \quad\quad
  \check{\bm \sigma}= \left\{ \begin{aligned}
     & {\bm \sigma} &\quad {\text on} ~ \mathcal{F}_{K^{\ast}}^{t_n}, \\
     & {\bf 0} &\quad {\text on} ~ \mathcal{F}_{K^{\ast}}^{t_{n-1}}, \\
     & \bm \sigma + \alpha v A^{\delta}{\bf n}^x_{\Omega} &\quad {\text on} ~ \mathcal{F}_{K^{\ast},n}^D. \\
  \end{aligned} \right.
                          \end{equation}

Using (\ref{localvari}) and the defined fluxes, the variational problem of (\ref{localnonhomomacro}) is to: Find
$(v_{K^{\ast}}^{(1)},{\bm \sigma}_{K^{\ast}}^{(1)})\in {\bf V}^{(1)}_{K^{\ast}}$ such that
\be \label{localdg}
\mathcal{A}^{(1)}(v_{K^{\ast}}^{(1)},{\bm \sigma}_{K^{\ast}}^{(1)};w,{\bm \tau}) =
\ell^{(1)}(w,{\bm \tau})  \quad \forall (w,{\bm \tau}) \in {\bf V}^{(1)}_{K^{\ast}},
\en
where
\begin{eqnarray}
& \mathcal{A}^{(1)}(v,{\bm \sigma};w,{\bm \tau})=
-\int_{K^{\ast}} \bigg({\bm \sigma}\cdot ( \frac{\partial{\bm \tau}}{\partial t } + A^{\frac{1}{2}}\nabla w ) + v(\nabla\cdot A^{\frac{1}{2}}{\bm \tau} + c^{-2}\frac{\partial w}{\partial t} ) \bigg)~dV
\cr &
+  \int_{\mathcal{F}_{K^{\ast}}^{t_n}} \big( c^{-2}v w + {\bm \sigma} \cdot {\bm \tau}\big)~d{\bf x} + \int_{\mathcal{F}_{K^{\ast},n}^D} \bigg({\bm \sigma} \cdot w A^{\frac{1}{2}}{\bf n}^x_{K^{\ast}}  + \alpha v w A^{\delta}{\bf n}^x_{K^{\ast}} \cdot A^{\frac{1}{2}}{\bf n}^x_{K^{\ast}} \bigg) ~dS,
 \end{eqnarray}
 and \be \nonumber
 \ell^{(1)}(w,{\bm \tau})= \int_{K^{\ast}} f w ~dV.
 \en



 \subsection{Discretization of the variational problems}
We decompose the discrete solution $(v_h, ~{\bm \sigma}_h)$
of the problem \((\ref{model})\) into $(v_h, ~ {\bm \sigma}_h)= (v_h^{(1)}, ~{\bm \sigma}_h^{(1)}) + (v_h^{(2)}, ~{\bm \sigma}_h^{(2)})$,
where \((v_h^{(1)}, ~{\bm \sigma}_h^{(1)})\) defined later is the discrete solution of continuous variational formulations (\ref{localdg}), and \((v_h^{(2)}, ~{\bm \sigma}_h^{(2)})\in {\bf V}_h({\cal T}_h)\) is the discrete solution of continuous Trefftz DG variational formulation (\ref{nontrefftzcontivaria}); namely, find \((v_h^{(2)}, ~{\bm \sigma}_h^{(2)})\in {\bf V}_h({\cal T}_h)\) such that
\be \label{nonhomotrefftzvaria}
\mathcal{A}(v_h^{(2)},{\bm \sigma}_h^{(2)};w,{\bm \tau}) =
\tilde{\ell}(w,{\bm \tau})- \mathcal{A}(v_h^{(1)},{\bm \sigma}_h^{(1)};w,{\bm \tau})  \quad \forall (w,{\bm \tau}) \in {\bf V}_h({\cal T}_h).
\en

  Let $\mathbb{Q}_q(K^{\ast})$ denote the set of polynomials of the same degree $q$ in each of the $n+1$ variables. Define ${\bf V}^{(1)}_{h,K^{\ast}}= \mathbb{Q}_q(K^{\ast})^{1+d}$ and ${\bf V}^{(1)}_{h,K}= \mathbb{Q}_q(K)^{1+d}$.

Then a discretized version of the continuous variational problem
(\ref{localdg}): Find $(v^{(1)}_{h,K^{\ast}}, ~{\bm \sigma}^{(1)}_{h,K^{\ast}})\in
{\bf V}^{(1)}_{h,K^{\ast}}$ such that
\be \label{localdisdg}
\mathcal{A}^{(1)}(v_{h,K^{\ast}}^{(1)},{\bm \sigma}_{h,K^{\ast}}^{(1)};w,{\bm \tau}) =
\ell^{(1)}(w,{\bm \tau})  \quad \forall (w,{\bm \tau})\in {\bf V}^{(1)}_{h,K^{\ast}}.
\en

Define $(v_h^{(1)}, ~{\bm \sigma}_h^{(1)})\in \prod_{K\in {\cal T}_h} {\bf V}^{(1)}_{h,K}$ by $(v_h^{(1)}, ~{\bm \sigma}_h^{(1)})|_K=(v^{(1)}_{h,K^{\ast}}, ~{\bm \sigma}^{(1)}_{h,K^{\ast}})|_K$.

Determine each local fictitious domain $K_{\bf x}^{\ast}$ by using the inverse transformation of (\ref{tran1}) acting on $\hat K_{\hat{\bf x}}^{\ast}$.
 A natural way is to choose $\hat K_{\hat{\bf x}}^{\ast}$ as the geometric sphere, e.g.
 the disc for the two-dimensional case and the sphere for the three-dimensional case, whose radius and center are denoted by $r_{\hat K_{\hat{\bf x}}}$ and $O_{\hat K_{\hat{\bf x}}}$, respectively. Notice that the center and the radius can be calculated easily. Then the variational problems (\ref{localdisdg}) can be solved easily by using the polar coordinate transformation for the calculation of the involved integrations. We would like to emphasize that the discrete problems (\ref{localdisdg}) are local and independent each other for $K\in  {\cal T}_h$, so they can be explicitly solved in parallel and the cost is small for low dimensional space case.

\begin{remark}\label{remarkcombound}
We would like to point out, there is no constraint for our proposed method for the nonhomogeneous model on the data of homogeneous Neumann boundary conditions from the original IBVP, owing to the fact that the analytic solution of the artificially constructed local IBVP (\ref{localnonhomomacro}) automatically satisfies the homogeneous Neumann boundary condition (Note that $\mathcal{F}_{K^{\ast},n}^N=\emptyset$), which is necessary in the theoretical error analysis of the DG method \cite{BMPS}.
\end{remark}

\begin{remark}
As stated in \cite{hy3}, if the nonhomogeneous local problem is defined on each element $K$, which is a non-smooth domain, then the analytic solution has only low regularity even if the analytic solution of the original problem defined on the global solution domain is smooth enough. Surprisingly, if we choose another alternative to define nonhomogeneous local problems on each {\it nonoverlapping} time slab $D_n (1\leq n \leq N)$, the combined DG scheme can also preserve the same orders of convergence as for the local smooth overlapping case (see Section \ref{overlapsec}).
\end{remark}

\section{$A$: piecewise-constant matrix} \label{piecesec}
In this section, we consider the model (\ref{model}) in which $A$ is a piecewise-constant positive definite matrix. By the derivation of section 3 and section 4, we still employ the variational formulations (\ref{trefftzcontivaria}) and (\ref{trefftzvaria}) to solve the continuous field $(v,{\bm \sigma})$ and its approximation $(v_{h},{\bm \sigma}_{h})$, respectively. Note that, for the piecewise constant model, the coordinate transformation (\ref{tran1}) can map $\Omega$ into disjoint subdomains, for example, see Figure 2, where $A|_{\Omega_1}=\left(
                  \begin{array}{cc}
                    \frac{3}{2} & -\frac{1}{2} \\
                    -\frac{1}{2}  & \frac{3}{2} \\
                  \end{array}
                \right), A|_{\Omega_2}=\left(
                  \begin{array}{cc}
                    2 & -1 \\
                    -1  & 2 \\
                  \end{array}
                \right), A|_{\Omega_3}=\left(
                  \begin{array}{cc}
                    4 & -2 \\
                    -2  & 4 \\
                  \end{array}
                \right), A|_{\Omega_4}=\left(
                  \begin{array}{cc}
                    5 & -3 \\
                    -3  & 5 \\
                  \end{array}
                \right).$
 It seems impossible to make the transformed mesh satisfy the shape regular and quasi uniform hypothesis.

\vskip -0.1in
\begin{figure}[H]
\begin{center}
\begin{tabular}{cc}
 \epsfxsize=0.5\textwidth\epsffile{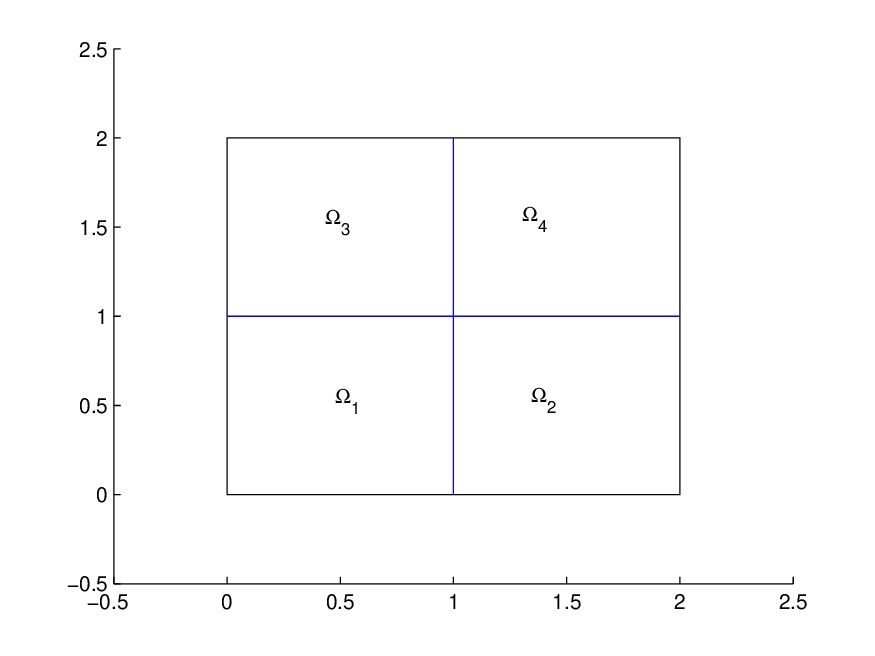}&
 \epsfxsize=0.5\textwidth\epsffile{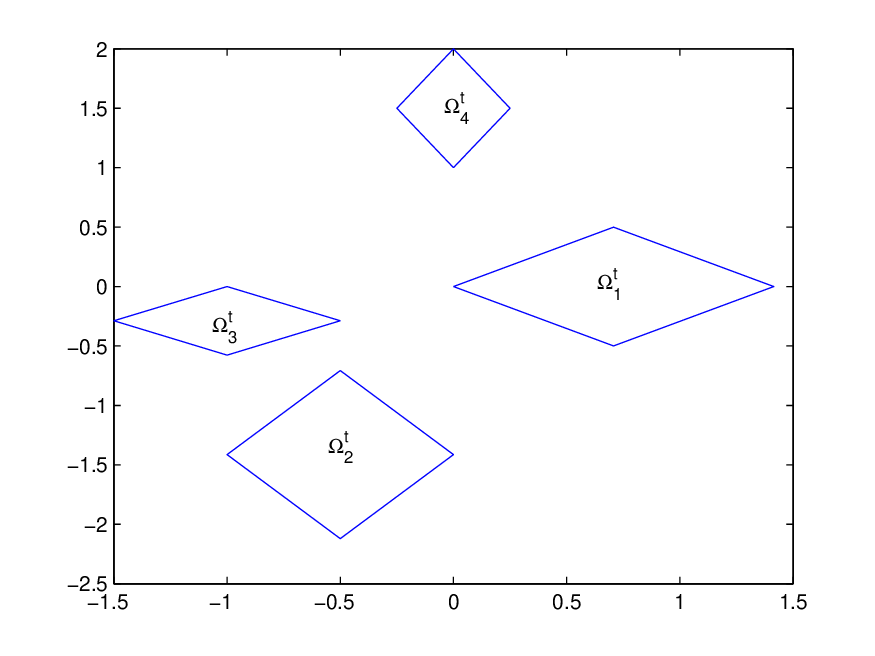}\\
\end{tabular}
\end{center}
\caption{A cube mapped into disjoint subdomains by the coordinate transformation (\ref{tran1}).}
\label{cube}
\end{figure}

 In order to make our algorithm adaptive to this model, we
divide $\Omega$ directly so that the mesh ${\cal T}_{h_{{\bf x}}}^{\bf x}=\{K_{\bf x}\}$ satisfies the shape regular and quasi uniform conditions. Naturally, compared with the partition introduced in section 5.1, by the simple and direct calculation, the proposed triangulation here satisfies that, for each subdomain where $A$ is a constant matrix,
\be \label{orihelmrelahhath}
c_0 ||\Lambda^{\frac{1}{2}}||^{-1}h_{{\bf x}} \leq \hat{h}_{\hat{\bf x}} \leq C_0||\Lambda^{-\frac{1}{2}}||h_{{\bf x}}, ~~\text{and} ~~c_0 ||\Lambda^{\frac{1}{2}}||^{-1}h \leq \hat{h} \leq C_0||\Lambda^{-\frac{1}{2}}||h.
\en

Of course, Lemma \ref{wellposedness} stating the existence and uniqueness of the Trefftz DG solution of (\ref{trefftzvaria}) and the upper bounds of the bilinear form (\ref{adg1}) still holds for this situation. In the following, we give the error estimates of Trefftz discontinuous Galerkin approximations generated by (\ref{trefftzvaria}).

\begin{theorem} \label{piecewise_anisoerr}
Assume that the IBVP solution $(v,{\bm \sigma})\in C^{k_t-1}(I;H^{k_{\bf x}+1}(\Omega))\times C^{k_t}(I;H^{k_{\bf x}}(\Omega)^d)$, and that $s=\text{min}\{p,k_t-1,k_{\bf x}-1\}$. Then we have,
\be\nonumber
\frac{1}{2} \bigg( ||c^{-1}  (v-v_h)||_{L^2(\Omega\times\{T\})} + ||\bm\sigma-\bm\sigma_h||_{L^2(\Omega\times\{T\})^d} \bigg)
 \leq  |||(v,{\bm \sigma}) - (v_h,{\bm \sigma}_{h})|||_{\text{DG}(Q)}
\leq  C\rho_{\text{max}}^{\frac{s+1}{2}}h^{s+\frac{1}{2}}
 |(v,{\bm \sigma})|_{H^{s+1}(Q)^{1+d}},
\en
where $\rho_{\text{max}}$ is the maximum of the condition number of piecewise constant matrices $A$.
\end{theorem}

{\it Proof}. By (\ref{abstarct}) and (\ref{transvari}), we obtain
\be \label{pieceerr1}
|||(v,{\bm \sigma}) - (v_h,{\bm \sigma}_h)|||_{\text{DG}(Q)} \leq
|||(v,{\bm \sigma}) - Q_h(v,\bm\sigma)|||_{\text{DG}(Q)^+}.
\en
By (\ref{errderi8}) and (\ref{errderi9}) posed on the triangulation ${\cal T}_h(Q)$, the scaling argument, the approximation estimate (\ref{isotreeappr4}), and (\ref{orihelmrelahhath}), it yields that
 \beq \label{pieceerr6} 
 &|||(v,{\bm \sigma}) - (v_h,{\bm \sigma}_h)|||_{\text{DG}(Q)}
  \leq  \sum\limits_{K=K_{x}\times I_{n'}\in {\cal T}_{h}(Q)}\bigg[
h^{-\frac{1}{2}}||(v,{\bm \sigma})-Q_{h}(v,{\bm \sigma})||_{L^2(K)^{1+d}}
\cr  &
+  h^{\frac{1}{2}}\big|\big((v,{\bm \sigma})-Q_{h}(v,{\bm \sigma})\big)\big|_{H^1(K)^{1+d}}\bigg]
\leq C\rho_{\text{max}}^{\frac{s+1}{2}}~h^{s+{\frac{1}{2}}}
|(v,{\bm \sigma})|_{H^{s+1}(Q)^{1+d}}. \quad\Box
 \eq

\begin{remark}
We would like to point that, although the error estimates in the Theorem \ref{piecewise_anisoerr} has the same $h$-convergence order as the estimates in the Theorem \ref{anisoerr}, the convergence order with respect to the condition number in the Theorem \ref{piecewise_anisoerr} is clearly lower than that in the Theorem \ref{anisoerr}, owing to the fact that the transformed mesh can not satisfy the shape regular and quasi uniform hypothesis compared with homogeneous media. But, one strategy to improve the convergence order with respect to $\rho$ for the proposed method is to employ nonconforming meshes in the original domain, which will be investigated in the next article.
\end{remark}

\section{Numerical experiments}

In this Section, we apply the proposed methods to solve the wave propagation in anisotropic media, and we report numerical results to verify the efficiency of the method.

The wave speed is fixed at $c=1$. As described in Section 4, we choose the same number \(p\) of basis functions for every elements \(\Omega_k\), and consider the following choice of numerical fluxes for the proposed methods: the constant parameters $\alpha =\beta = 1$. Meanwhile, we compare numerical performances among the Trefftz methods and the high-order DG method introduced in \cite{BMPS}. The new proposed assumption on the shape regularity of polyhedral meshes $\hat{{\cal T}}_{\hat{h}_{\hat{\bf x}}}$ is employed for the high-order DG method. For the convenience of statement, we call the method (\ref{trefftzcontivaria}) as ``Method-I", and (\ref{isotrefftzcontivaria}) described in Section 6 as  ``Method-II".

In our tests, we estimate the convergence orders of the approximations by measuring the relative errors in $L^2(\Omega)$-norm at a given time $T$, and the errors in $|||\cdot|||_{\text{DG}}-$norm, respectively. All of the computations have been done in MATLAB, and the system matrix was computed by numerical integration. ``DOFs" represents the number of degree of freedoms equal to the elements multiplied by the number of basis functions per element.


\subsection{Homogeneous case}
We use uniform mesh with $ h_{{\bf x}} \approx h_t =2^{-l}, l\in \mathbb{N}$.
In the Tables \ref{homorhofixtoh}-\ref{3dhomohfixtorho}, the collum labelled ``Error" shows the numerical errors in  relative $L^2(\Omega\times\{T\})-$norm and $|||\cdot|||_{\text{DG}}-$norm: if $p=1$, it is given for the mesh level $l=4$, else if $p\in \{2,3\}$ for $l=3$.

\subsubsection{2D space case with Neumann boundary conditions}

We choose the space-time domain $Q=(0,1)^2\times(0,1)$, and set the anisotropic matrix $A = \left(
  \begin{array}{cc}
  \lambda_1 a^2+ \lambda_2b^2 & ab(\lambda_2-\lambda_1)  \\
  ab(\lambda_2-\lambda_1)  & \lambda_1b^2+ \lambda_2a^2\\
  \end{array}
\right)$, where $a=\frac{1}{\sqrt{2}},~b=\frac{1}{\sqrt{2}}$, and $0<\lambda_1<1,\lambda_2=1$.
 Consider the exact smooth solution
 \beq
\hat U(\hat{\bf x},t) & = & \text{sin}(\pi \hat x_1) \text{sin}(\pi \hat x_2)\text{sin}(\sqrt{2}\pi t),
\cr v &=& \hat v = \frac{\partial \hat U}{\partial t}, ~  {\bm \sigma}= P^T\hat{\bm \sigma} = - P^T\hat\nabla \hat U,
\eq
to the IBVP (\ref{model}), with nonhomogeneous Neumann boundary conditions.

The convergence rates with respect to $h$ are given in the Table \ref{homorhofixtoh}.
\begin{center}
       \tabcaption{}
\label{homorhofixtoh} \vskip -0.3in
       Convergence rates of the TDG schemes w.r.t. $h$.
\begin{tabular}{|c|c|cc|cc|cc|} \hline
   &   &  \multicolumn{2}{|c|}{$v_h$} &  \multicolumn{2}{|c|}{$\bm\sigma_h$} &  \multicolumn{2}{|c|}{$|||\cdot|||_{\text{DG}}$} \\ \hline
 \text{Method-I}  & $\rho$ & \text{Error} & \text{Rate} & \text{Error} & \text{Rate} &\text{Error} & \text{Rate} \\ \hline
     \multirow{3}*{\(p=1\)} & 2 & 4.70e-2  & 2.30  & 8.27e-3   & 2.41  & 2.02e-1 &  1.49 \\ \cline{2-8}
   & 4 &  4.31e-2 & 2.31  & 8.25e-3   &  2.39 &   1.93e-1 & 1.49  \\ \cline{2-8}
   & 16 & 4.60e-2  &  2.30 & 8.45e-3   & 2.33  &  1.94e-1 &  1.48  \\ \hline
    \multirow{3}*{\(p=2\)} & 2 & 1.16e-2 & 3.38  & 2.56e-3   & 3.17  &  5.49e-2 &   2.51 \\ \cline{2-8}
   & 4 & 1.04e-2  &  3.31 & 2.51e-3   & 3.16 &  5.03e-2 &  2.52 \\ \cline{2-8}
   & 16 & 1.18e-2  & 3.25   & 2.86e-3   & 3.14  & 5.75e-2  &  2.58\\ \hline
     \multirow{3}*{\(p=3\)} & 2 &  7.55e-4 & 4.12 & 2.51e-4   & 3.90  &  3.83e-3 &   3.59 \\ \cline{2-8}
   & 4 &  6.77e-4  &  4.10 & 2.39e-4   & 3.92   & 3.41e-3  &  3.59 \\ \cline{2-8}
   & 16 & 7.96e-4 &4.06  &  2.87e-4  & 3.96  & 3.33e-3 & 3.55 \\ \hline
   \text{Method-II}  & $\rho$ & \text{Error} & \text{Rate} & \text{Error} & \text{Rate} &\text{Error} & \text{Rate} \\ \hline
     \multirow{3}*{\(p=1\)} & 2  & 4.85e-2 & 2.30  & 8.68e-3   & 2.43  &  2.01e-1 & 1.49 \\ \cline{2-8}
   & 4 & 4.52e-2 & 2.30  &  8.88e-3   & 2.43  &  1.91e-1 &  1.49 \\ \cline{2-8}
   & 16 & 5.02e-2  & 2.31 & 9.24e-3  &  2.41 & 1.91e-1 & 1.48 \\ \hline
    \multirow{3}*{\(p=2\)} & 2 & 1.20e-2 & 3.44 & 2.58e-3    & 3.18  & 5.47e-2  & 2.52\\ \cline{2-8}
   & 4 & 1.10e-2 &  3.40 & 2.53e-3   & 3.18  & 4.99e-2  & 2.52 \\ \cline{2-8}
   & 16 & 1.27e-2 & 3.43  & 2.83e-3   & 3.21  &  5.17e-2 & 2.51 \\ \hline
     \multirow{3}*{\(p=3\)} & 2 & 7.69e-4  & 4.13  & 2.51e-4    &  3.89  &  3.82e-3 & 3.49 \\ \cline{2-8}
   & 4 & 7.02e-4  & 4.13  &  2.38e-4  & 3.89  & 3.39e-3 &  3.50 \\ \cline{2-8}
   & 16 & 8.37e-4  & 4.11  & 2.84e-4    & 3.89  & 3.59e-3  &  3.49 \\ \hline
   \end{tabular}
     \end{center}


We can obtain that $|| v-v_h||_{L^2(\Omega\times T)} \approx ||\bm\sigma-\bm\sigma_h||_{L^2(\Omega\times T)^2} =\mathcal{O}(h^{p+1})$. The last column showing the experimental convergence rates of the errors measured in $|||\cdot|||_{\text{DG}}-$norm indicates that the estimates of Theorems \ref{anisoerr} and \ref{isoerr} are sharp.

The convergence rates with respect to $\rho$ are given in the Table \ref{homohfixtorho}. 


\begin{center}
       \tabcaption{}
\label{homohfixtorho} \vskip -0.3in
       Convergence rates of the TDG schemes w.r.t. $\rho$.
\begin{tabular}{|c|c|c|c|c|c|c|c|} \hline
     & & \multicolumn{2}{|c|}{$v_h$} &  \multicolumn{2}{|c|}{$\bm\sigma_h$} &  \multicolumn{2}{|c|}{$|||\cdot|||_{\text{DG}}$} \\ \hline
   \text{Method-I} & $\rho$ &   \text{Error} & \text{Rate} & \text{Error} & \text{Rate} &\text{Error} & \text{Rate} \\ \hline
 \multirow{3}*{$p=1$} &  $32$ & 4.61e-2  &   & 8.12e-3  &     &  1.92e-1 &  \\ \cline{2-8} 
  & 64  & 4.84e-2 &  0.0702 &  8.53e-3  &  0.0711 & 1.96e-1 &  0.0297  \\ \cline{2-8}
   & 128  & 5.01e-2  & 0.0498 & 8.89e-3   &  0.0596  & 1.97e-1  & 0.0073   \\ \hline
 \multirow{3}*{$p=2$} &  $32$ & 1.15e-2 &   & 2.94e-3  &     & 5.16e-2  &  \\ \cline{2-8} 
  & 64  & 1.28e-2 & 0.1545 & 3.35e-3   &  0.1883 & 5.43e-2 &  0.0736  \\ \cline{2-8}
   & 128 & 1.36e-2  &  0.0875  & 3.77e-3   &  0.1704
  & 5.62e-2 &  0.0496  \\ \hline
    \multirow{3}*{$p=3$} & 8 & 6.87e-4 &  & 2.54e-4   &   &  3.37e-3 &    \\ \cline{2-8}
   & 16 & 7.96e-4 & 0.2125 & 2.87e-4   & 0.1762  & 3.61e-3 & 0.0993 \\ \cline{2-8}
      & $32$ & 7.89e-4 &  -0.0127 & 2.96e-4  &  0.0445   &  3.54e-3 & -0.0282 
    \\ \hline
     \text{Method-II} & $\rho$ &   \text{Error} & \text{Rate} & \text{Error} & \text{Rate} &\text{Error} & \text{Rate} \\ \hline
 \multirow{3}*{$p=1$} &  $32$ &  5.06e-2 &   & 8.95e-3  &     & 1.89e-1  &  \\ \cline{2-8} 
  & 64  & 5.35e-2 & 0.0804 & 9.52e-3   & 0.0891  & 1.92e-1  &    0.0227 \\ \cline{2-8}
   & 128 & 5.58e-2  & 0.0607 &  9.99e-3  &  0.0695  & 1.94e-1 &  0.0150  \\ \hline
 \multirow{3}*{$p=2$} &  $32$ & 1.25e-2 &   & 2.88e-3  &     & 5.08e-2  &  \\ \cline{2-8} 
  & 64  & 1.40e-2  & 0.1635  & 3.06e-3   &  0.0875  & 5.34e-2 & 0.0720   \\ \cline{2-8}
   & 128  & 1.49e-2 & 0.0899  & 3.21e-3   & 0.0690  & 5.52e-2  & 0.0478  \\ \hline
    \multirow{3}*{$p=3$} &  8 & 7.20e-4 &  &  2.52e-4  &   & 3.35e-3  &     \\ \cline{2-8} 
  & 16 & 8.37e-4 &  0.2172 &2.84e-4     &  0.1725 & 3.59e-3 &   0.0998  \\ \cline{2-8}
   &  $32$ & 8.36e-4  &  -0.0017 & 2.94e-4  &   0.0499   & 3.52e-3  &  -0.0284
  \\ \hline
   \end{tabular}
     \end{center}

We can obtain that $|| v-v_h||_{L^2(\Omega\times T)}, ||\bm\sigma-\bm\sigma_h||_{L^2(\Omega\times T)^2}, |||(v,{\bm \sigma}) - ( v_h,{\bm \sigma}_{h})|||_{\text{DG}} \lesssim C(v,\bm\sigma)\mathcal{O}(\rho^{\frac{1}{4}})$, which indicates that the estimates of Theorems \ref{anisoerr} and \ref{isoerr} are sharp.

\subsubsection{3D space case with Neumann boundary conditions}
We choose the space-time domain $Q=(0,1)^3\times(0,1)$, and set the anisotropic matrix $A = \left(
  \begin{array}{ccc}
  \lambda_1a^2+ \lambda_2b^2 & ab(\lambda_2-\lambda_1) & 0  \\
  ab(\lambda_2-\lambda_1)  & \lambda_1b^2+ \lambda_2a^2 & 0\\
  0 & 0 & 1 \\
  \end{array}
\right)$, where $a=\frac{1}{\sqrt{2}},~b=\frac{1}{\sqrt{2}}$, and $0<\lambda_1,\lambda_2<1$.
Consider the exact smooth solution
 \beq
\hat U(\hat{\bf x},t) & = & \text{sin}(\pi \hat x_1) \text{sin}(\pi \hat x_2) \text{sin}(\pi \hat x_3)\text{sin}(\sqrt{3}\pi t),
\cr v &=& \hat v = \frac{\partial \hat U}{\partial t}, ~  {\bm \sigma}= P^T\hat{\bm \sigma} = - P^T\hat\nabla \hat U,
\eq
to the IBVP (\ref{model}), with nonhomogeneous Neumann boundary conditions.


The convergence rates with respect to $h$ are given in the Table \ref{new3dhomorhofixtoh}. The convergence rates with respect to $\rho$ are given in the Table \ref{3dhomohfixtorho}.



\begin{center}
       \tabcaption{}
\label{new3dhomorhofixtoh} \vskip -0.3in
       Convergence rates of the TDG schemes w.r.t. $h$.
\begin{tabular}{|c|c|cc|cc|cc|} \hline
   &   &  \multicolumn{2}{|c|}{$v_h$} &  \multicolumn{2}{|c|}{$\bm\sigma_h$} &  \multicolumn{2}{|c|}{$|||\cdot|||_{\text{DG}}$} \\ \hline
 ($p,\rho$)&\text{Method} & \text{Error} & \text{Rate} & \text{Error} & \text{Rate} &\text{Error} & \text{Rate} \\ \hline
    \multirow{2}*{(1,2)} & \text{-I} &  2.41e-2   & 2.35 &  7.19e-3  &  2.53 & 1.68e-1  & 1.57 \\ \cline{2-8}
   & \text{-II} & 2.39e-2  &  2.37 &  7.45e-3  &  2.54 &  1.62e-1 &  1.59   \\ \hline
     \multirow{2}*{(2,2)} & \text{-I} &1.03e-2 & 3.41 &   7.54e-3 & 3.26  & 8.16e-2  &  2.51 \\ \cline{2-8}
   & \text{-II} & 1.08e-2 & 3.41  &  7.78e-3  & 3.25  &  8.02e-2 &  2.53  \\ \hline
     \multirow{2}*{(3,2)} & \text{-I} & 7.76e-4 & 4.12 & 8.45e-4   & 4.05  & 8.62e-3  &  3.52 \\ \cline{2-8}
   & \text{-II} & 7.94e-4  & 4.11 &  8.75e-4  & 4.02  & 8.85e-3  &  3.48  \\ \hline
   \end{tabular}
     \end{center}

\begin{center}
       \tabcaption{}
\label{3dhomohfixtorho}  \vskip -0.3in
      Convergence rates of the TDG schemes w.r.t. $\rho$.
\begin{tabular}{|c|c|c|c|c|c|c|c|} \hline
     & & \multicolumn{2}{|c|}{$v_h$} &  \multicolumn{2}{|c|}{$\bm\sigma_h$} &  \multicolumn{2}{|c|}{$|||\cdot|||_{\text{DG}}$} \\ \hline
  \text{Method-I} & $\rho$ &   \text{Error} & \text{Rate} & \text{Error} & \text{Rate} &\text{Error} & \text{Rate} \\ \hline
 \multirow{3}*{$p=1$} &  $4$ & 1.53e-2  &   & 7.83e-3  &     & 1.05e-1  &  \\ \cline{2-8} 
  & 8  & 1.48e-2 &  -0.0479  &  8.02e-3  &  0.0346 & 1.06e-1 &  0.0137   \\ \cline{2-8}
   & 16 & 1.48e-2  & 0 &  8.58e-3  &  0.0974 & 1.09e-1 &  0.0403  \\ \hline
 \multirow{3}*{$p=2$} &  $4$ & 9.67e-3 &   &  4.97e-3 &     & 4.63e-2  &  \\ \cline{2-8} 
  & 8  & 9.48e-3  &  -0.0286 & 4.78e-3   &  -0.0562  & 4.63e-2 &  0  \\ \cline{2-8}
   & 16  & 9.73e-3  & 0.0376 &  4.88e-3  &   0.0299 & 4.59e-2  & -0.0125   \\ \hline
    \multirow{3}*{$p=3$} &  $4$ & 9.18e-4  &   & 7.07e-4  &     & 5.34e-3  &  \\ \cline{2-8} 
  & 8 & 8.98e-4  & -0.0318 & 7.03e-4   &  -0.0082 & 5.34e-3  & 0   \\ \cline{2-8}
   & 16 & 9.57e-4  & 0.0918 & 7.46e-4   & 0.0857  & 5.45e-3  &  0.0294  \\ \hline
 \text{Method-II} & $\rho$ &   \text{Error} & \text{Rate} & \text{Error} & \text{Rate} &\text{Error} & \text{Rate} \\ \hline
 \multirow{3}*{$p=1$} &  $4$ & 1.52e-2 &   & 8.14e-3  &     & 1.03e-1  &  \\ \cline{2-8} 
  & 8  & 1.45e-2  &-0.0680  & 8.67e-3   &  0.0910 & 1.03e-1  & 0   \\ \cline{2-8}
   & 16 & 1.44e-2 & -0.0100 &  9.66e-3  &   0.1560 & 1.01e-1 &   -0.0283  \\ \hline
 \multirow{3}*{$p=2$} &  $4$ &  1.04e-2 &   & 4.86e-3  &     &  4.56e-2 &  \\ \cline{2-8} 
  & 8  & 1.06e-2 &  0.0275 &   4.66e-3 &  -0.0606  & 4.59e-2  &  0.0095  \\ \cline{2-8}
   & 16  & 1.13e-2 & 0.0923 &  4.75e-3  & 0.0276  & 4.59e-2 &  0  \\ \hline
    \multirow{3}*{$p=3$} &  $4$ & 9.92e-4 &   & 7.07e-4  &     & 5.34e-3  &  \\ \cline{2-8} 
  & 8 & 1.01e-3 & 0.0259 &  7.07e-4   &  0 & 5.36e-3  &  0.0054  \\ \cline{2-8}
   & 16 & 1.10e-3  &  0.1231 &  7.42e-4  &  0.0697 & 5.51e-3  &  0.0398  \\ \hline
      \end{tabular}
     \end{center}

We can obtain that $|| v-v_h||_{L^2(\Omega\times T)} \approx ||\bm\sigma-\bm\sigma_h||_{L^2(\Omega\times T)^2} =\mathcal{O}(h^{p+1})$, and $|| v-v_h||_{L^2(\Omega\times T)}, ||\bm\sigma-\bm\sigma_h||_{L^2(\Omega\times T)^2}, |||(v,{\bm \sigma}) - ( v_h,{\bm \sigma}_{h})|||_{\text{DG}} \lesssim C(v,\bm\sigma)\mathcal{O}(\rho^{\frac{1}{4}})$, which indicate that the estimates of Theorems \ref{anisoerr} and \ref{isoerr} are sharp.


\subsection{Nonhomogeneous case}
We use uniform mesh with $ h_{{\bf x}} \approx h_t =2^{-l}, l\in \mathbb{N}$.

\subsubsection{1D space case for smooth solution with Dirichlet boundary conditions}

We choose the space-time domain $Q=(0,1)\times(0,1)$. Consider the exact smooth solution
 \beq
U(x,t) & = & \text{sin}(\pi x) \text{sin}(\sqrt{2}\pi t),
\cr v& =& \partial_t U, ~\sigma = -\partial_x U,
\eq
to the IBVP (\ref{model}), with $A=1$, nonhomogeneous soure $f$:
$$f=-\pi^2\text{sin}(\pi x) \text{sin}(\sqrt{2}\pi t),$$
and homogeneous Dirichlet boundary conditions.

The convergence rates are given in the Table \ref{nnt1} for different choices of $p$ and $q$. 


\begin{center}
       \tabcaption{}
\label{nnt1} \vskip -0.3in
       Convergence rates of the combined space-time DG scheme w.r.t. $h$. 
\begin{tabular}{|c|c|c|c|c|c|c|c|} \hline
   &   &  \multicolumn{2}{|c|}{$v_h$} &  \multicolumn{2}{|c|}{$\sigma_h$} &  \multicolumn{2}{|c|}{$|||\cdot|||_{\text{DG}}$} \\ \hline
  ($p,q$) & $h$ & \text{Error} & \text{Rate} & \text{Error} & \text{Rate} &\text{Error} & \text{Rate} \\ \hline
     \multirow{3}*{\((1,1)\)} & 1/8 & 5.32e-2  &  & 4.63e-2  &    & 3.64e-1  &      \\ \cline{2-8}
   & 1/16 & 1.22e-2 & 2.24  & 8.70e-3  & 2.24   & 1.31e-1  &  1.47   \\ \cline{2-8}
   & 1/32 & 2.98e-3  & 2.11 &  1.77e-3 &  2.11  & 4.66e-2  &   1.49   \\ \hline
      \multirow{3}*{\((1,2)\)} & 1/8 & 5.19e-2 &  & 4.78e-2  &    &  3.30e-1 &       \\ \cline{2-8}
   & 1/16 & 1.21e-2  & 2.22 & 9.08e-3  &  2.22  & 1.21e-1   &   1.45    \\ \cline{2-8}
   & 1/32 & 2.98e-3 & 2.11 &  1.83e-3 &  2.11  & 4.31e-2  &  1.49   \\ \hline
     \multirow{3}*{\((2,1)\)} & 1/8 & 9.55e-3 &  & 9.11e-3  &    & 1.51e-1  &     \\ \cline{2-8}
   & 1/16 & 2.10e-3 & 2.33 & 1.58e-3  & 2.33   & 5.36e-2   &  1.50    \\ \cline{2-8}
   & 1/32 & 5.04e-4  & 2.15  &  3.10e-4  & 2.15   &   1.90e-2 &  1.50   \\ \hline
     \multirow{3}*{\((2,2)\)} & 1/4 & 2.70e-2 &  &  2.95e-2 &    & 1.07e-1  &       \\ \cline{2-8}
   & 1/8 & 3.02e-3 & 3.20 & 3.17e-3  &  3.20     &  1.99e-2 &  2.52   \\ \cline{2-8}
   & 1/16 & 3.50e-4 & 3.13 &  3.62e-4 &  3.13     & 3.57e-3  &  2.48  \\ \hline
      \multirow{3}*{\((2,3)\)}   & 1/4 & 2.72e-2  &  & 3.02e-2  &       & 9.52e-2  &   \\ \cline{2-8}
      & 1/8 & 3.04e-3 & 3.21 & 3.24e-3  &  3.21  &  1.73e-2 &  2.46     \\ \cline{2-8}
   & 1/16 &  3.55e-4 & 3.12 &  3.72e-4 &   3.12 & 3.72e-3  &  2.52     \\ \hline
      \multirow{3}*{\((3,2)\)} & 1/4 & 4.57e-3 &  & 8.12e-3  &    &  6.10e-2 &       \\ \cline{2-8}
   & 1/8 & 6.26e-4 & 3.21 & 7.99e-4  &  3.22  & 1.08e-2  &    2.50   \\ \cline{2-8}
   & 1/16 & 8.45e-5 & 2.99 &9.61e-5   & 2.99   & 1.91   &  2.50  \\ \hline
   \multirow{3}*{\((3,3)\)} & 1/2 & 5.90e-2 &  & 3.59e-2  &    & 1.10e-1  &       \\ \cline{2-8}
   & 1/4 & 2.69e-3 & 4.06  &  3.14e-3 & 4.04    &  9.02e-3 &  3.61     \\ \cline{2-8}
   & 1/8 & 1.29e-4 & 4.28  & 1.68e-4  &  4.29  & 7.69e-4  &  3.55    \\ \hline
   \end{tabular}
     \end{center}


From the Table \ref{nnt1}, we obtain that, as typical for DG methods, in the case of a regular enough solution, we observe the convergence orders of the errors in $L^2(\Omega\times\{T\})-$norm with the rate $\mathcal{O}(h^{\text{min}\{p+1,q+1\}})$,
and that the convergence rates of errors in $|||\cdot|||_{\text{DG}}-$norm are $O(h^{\text{min}\{p+\frac{1}{2},q+\frac{1}{2}\}})$, which support convergence rate optimality on the uniform refined meshes of the combined numerical DG scheme just as for the homogeneous case. Furthermore, the Table \ref{nnt1} shows that for the $L^2(\Omega\times\{T\})$ norms and  $|||\cdot|||_{\text{DG}}$ norms of errors generated by the proposed method, the choice $p=q+1$ is preferable compared with the choice $p=q$.

\subsubsection{2D space case for smooth solution with Neumann boundary conditions}
\label{test2d1t}
We choose the space-time domain $Q=(0,1)^2\times(0,1)$, and set the anisotropic matrix $A = \left(
  \begin{array}{cc}
  \lambda_1a^2+ \lambda_2b^2 & ab(\lambda_2-\lambda_1)  \\
  ab(\lambda_2-\lambda_1)  & \lambda_1b^2+ \lambda_2a^2\\
  \end{array}
\right)$, where $a=\frac{1}{\sqrt{2}},~b=\frac{1}{\sqrt{2}}$, and $0<\lambda_1<1,\lambda_2=1$.
 Consider the exact smooth solution
 \beq
U({\bf x},t) & = & \text{sin}(\pi x_1) \text{sin}(\pi x_2)\text{sin}(\sqrt{3}\pi t),
\cr v &=& \frac{\partial U}{\partial t}, ~  {\bm \sigma}=-A^{\frac{1}{2}}\nabla U,
\eq
to the IBVP (\ref{model}), with nonhomogeneous soure $f$ and nonhomogeneous Neumann boundary conditions.


In the Tables \ref{nnt2d}-\ref{hfixtorho}, the collum labelled ``Error" shows the numerical errors in  relative $L^2(\Omega\times\{T\})-$norm and $|||\cdot|||_{\text{DG}}-$norm: if $p=\in \{1,2\}$, it is given for the mesh level $l=4$, else if $p\in \{3,4\}$ for $l=3$. The convergence rates with respect to $h$ are given in the Table \ref{nnt2d} for different choices of $p=q+1$. 


\begin{center}
       \tabcaption{}
\label{nnt2d} \vskip -0.3in
       Convergence rates of the combined space-time DG scheme w.r.t. $h$.
\begin{tabular}{|c|cc|cc|cc|} \hline
      &  \multicolumn{2}{|c|}{$v_h$} &  \multicolumn{2}{|c|}{$\bm\sigma_h$} &  \multicolumn{2}{|c|}{$|||\cdot|||_{\text{DG}}$} \\ \hline
 $(p,q)$  & \text{Error} & \text{Rate} & \text{Error} & \text{Rate} &\text{Error} & \text{Rate} \\ \hline
    \((2,1)\) &  1.38e-2  & 3.32  & 1.23e-2   &  3.13 & 1.56e-1 &  1.87 \\ \hline
      \((3,2)\) & 1.43e-3 & 3.95  & 2.11e-3   &  3.42 & 2.22e-2 & 2.53 \\\hline
     \((4,3)\) & 5.55e-5  & 4.80  &  6.24e-5  & 5.07  & 8.49e-4 & 3.67   \\\hline
   \end{tabular}
     \end{center}


From the Table \ref{nnt2d}, we observe the convergence orders of  the errors 
in $L^2(\Omega\times\{T\})-$norm between $\mathcal{O}(h^{p+\frac{1}{2}})$ and $\mathcal{O}(h^{p+1})$. Besides, the convergence rates of errors in $|||\cdot|||_{\text{DG}}-$norm are $O(h^{\text{min}\{p+\frac{1}{2},q+\frac{1}{2}\}})$, which support convergence rate optimality on the uniform refined meshes of the combined numerical DG scheme just as for the homogeneous case.


The convergence rates with respect to $\rho$ are given in the Table \ref{hfixtorho}.


\begin{center}
       \tabcaption{}
\label{hfixtorho} \vskip -0.3in
      Convergence rates of the TDG schemes w.r.t. $\rho$.
\begin{tabular}{|c|c|c|c|c|c|c|c|} \hline
 &    &  \multicolumn{2}{|c|}{$v_h$} &  \multicolumn{2}{|c|}{$\bm\sigma_h$} &  \multicolumn{2}{|c|}{$|||\cdot|||_{\text{DG}}$} \\ \hline
($p,q$)& $\rho$   &  \text{Error} & \text{Rate} & \text{Error} & \text{Rate} &\text{Error} & \text{Rate} \\ \hline
 \multirow{3}*{$(1,0)$} &  $16$ & 3.59e-2  &   & 7.42e-2  &     &  6.58e-1   &  \\ \cline{2-8} %
  & 32  & 5.55e-2 &  0.6285 & 8.32e-2& 0.1652  & 5.68e-1  &   -0.2122  \\ \cline{2-8}
   & 64 &  5.51e-2 & -0.0104  &  9.87e-2  &   0.2465 & 5.30e-1 &   -0.0999  \\ \hline
\multirow{3}*{$(2,1)$} &  $16$ &  2.63e-2 &   & 3.40e-2  &     & 3.86e-1  &  \\ \cline{2-8} %
  & 32  & 2.93e-2 &0.1558  &  4.01e-2  & 0.2381  &  4.01e-1 &  0.0550  \\ \cline{2-8}
   & 64 &3.27e-2 & 0.1584 &  4.68e-2  &  0.2229  & 4.24e-1 &  0.0805  \\ \hline
\multirow{3}*{$(3,2)$} &  $2$ &  3.25e-3   &  & 8.01e-3  & &  9.81e-2 &  \\ \cline{2-8} %
  & 4  & 3.77e-3 & 0.2141  &  8.97e-3  &  0.1633  & 1.04e-1  &   0.0843 \\ \cline{2-8}
   & 8 & 3.84e-3  & 0.0265  &  9.26e-3  &   0.0459 & 1.05e-1  &   0.0138   \\ \hline
   \end{tabular}
     \end{center}


We can see that $|| v-v_h||_{L^2(\Omega\times T)}, ||\bm\sigma-\bm\sigma_h||_{L^2(\Omega\times T)^2}, |||(v,{\bm \sigma}) - ( v_h,{\bm \sigma}_{h})|||_{\text{DG}} \lesssim C(v,\bm\sigma)\mathcal{O}(\rho^{\frac{1}{4}})$, which coincides with the optimal convergence rates of the errors with respect to $\rho$ indicated by the Theorem \ref{anisoerr} for the homogeneous case.

\subsubsection{3D space case for smooth solution with Neumann boundary conditions}
We choose the space-time domain $Q=(0,1)^3\times(0,1)$, and set the anisotropic matrix $A = \left(
  \begin{array}{ccc}
  \lambda_1a^2+ \lambda_2b^2 & ab(\lambda_2-\lambda_1) & 0  \\
  ab(\lambda_2-\lambda_1)  & \lambda_1b^2+ \lambda_2a^2 & 0\\
  0 & 0 & 1 \\
  \end{array}
\right)$, where $a=\frac{1}{\sqrt{2}},~b=\frac{1}{\sqrt{2}}$, and $0<\lambda_1,\lambda_2<1$.
Consider the exact smooth solution
 \beq
U({\bf x},t) & = & \text{sin}(\pi x_1) \text{sin}(\pi x_2) \text{sin}(\pi x_3)\text{sin}(2\pi t),
\cr v& =& \frac{\partial U}{\partial t}, ~  {\bm \sigma}=-A^{\frac{1}{2}}\nabla U,
\eq
to the IBVP (\ref{model}), with nonhomogeneous soure $f$ and nonhomogeneous Neumann boundary conditions.

The convergence rates with respect to $h$ are given in the Table \ref{nnt3d} for different choices of $p=q+1$. The collum labelled ``Error" shows the numerical errors in  relative $L^2(\Omega\times\{T\})-$norm and $|||\cdot|||_{\text{DG}}-$norm: if $p=1$, it is given for the mesh level $l=4$, else if $p\in \{2,3\}$ for $l=3$.


\begin{center}
       \tabcaption{}
\label{nnt3d} \vskip -0.3in
       Convergence rates of the combined space-time DG scheme w.r.t. $h$.
\begin{tabular}{|c|cc|cc|cc|} \hline
   &     \multicolumn{2}{|c|}{$v_h$} &  \multicolumn{2}{|c|}{$\bm\sigma_h$} &  \multicolumn{2}{|c|}{$|||\cdot|||_{\text{DG}}$} \\ \hline
  ($p,q$)  & \text{Error} & \text{Rate} & \text{Error} & \text{Rate} &\text{Error} & \text{Rate} \\ \hline
      \((1,0)\) &
      6.96e-2  & 1.85  &  1.20e-1  & 1.21  & 2.07e-1 & 1.51     \\\hline
    \((2,1)\) &  2.10e-2 & 2.71  &  5.23e-2  & 2.07  & 2.22e-1  &  1.47   \\\hline
  \((3,2)\) & 2.97e-3   & 3.53  &  6.26e-3  &  3.17  & 2.87e-2 & 2.55    \\\hline
   \end{tabular}
     \end{center}

In the Table \ref{nnt3d}, we can see that, the convergence orders of  the errors 
in $L^2(\Omega\times\{T\})-$norm are between $\mathcal{O}(h^{p})$ and $\mathcal{O}(h^{p+1})$. Besides, the convergence rates of errors in $|||\cdot|||_{\text{DG}}-$norm are $O(h^{\text{min}\{p+\frac{1}{2},q+\frac{1}{2}\}})$, which support convergence rate optimality on the uniform refined meshes of the combined numerical DG scheme just as for the homogeneous case.


The convergence rates with respect to $\rho$ are given in the Table \ref{3d1thfixtorho}.


\begin{center}
       \tabcaption{}
\label{3d1thfixtorho}
\vskip -0.3in
      Convergence rates of the TDG schemes w.r.t. $\rho$.
\begin{tabular}{|c|c|c|c|c|c|c|c|} \hline
 &    &  \multicolumn{2}{|c|}{$v_h$} &  \multicolumn{2}{|c|}{$\bm\sigma_h$} &  \multicolumn{2}{|c|}{$|||\cdot|||_{\text{DG}}$} \\ \hline
($p,q$)& $\rho$   &  \text{Error} & \text{Rate} & \text{Error} & \text{Rate} &\text{Error} & \text{Rate} \\ \hline
 \multirow{3}*{$(1,0)$} &  $4$ & 7.47e-2 &  &  7.41e-2  &   & 1.13e-1 &  \\ \cline{2-8} %
  & 8  & 7.25e-2 &  -0.0431 &  7.86e-2  &  0.0851 & 1.09e-1 & -0.0520   \\ \cline{2-8}
   & 16 & 7.07e-2 &  -0.0362 &  8.11e-2  &  0.0452  & 1.06e-1 &  -0.0403  \\ \hline
\multirow{3}*{$(2,1)$} &  $4$ & 2.39e-2  &  & 5.69e-2   &   & 2.54e-1 &   \\ \cline{2-8} %
  & 8  &2.66e-2  & 0.1544 &  6.56e-2  &  0.2053 & 2.72e-1 & 0.0988  \\ \cline{2-8}
   & 16 & 3.01e-2 & 0.1783 &  7.73e-2  & 0.2368  & 2.98e-1 & 0.1317   \\ \hline
\multirow{3}*{$(3,2)$} &  $4$ & 3.03e-3  &   & 7.24e-3  &     &  3.01e-2 &  \\ \cline{2-8} %
  & 8  & 3.49e-2 & 0.2039 &  8.25e-3  & 0.1884  & 3.34e-2 & 0.1501   \\ \cline{2-8}
   & 16 & 4.05e-2 & 0.2147  & 9.42e-3   & 0.1913  & 3.69e-2 & 0.1438   \\ \hline
   \end{tabular}
     \end{center}


We can see that $|| v-v_h||_{L^2(\Omega\times T)}, ||\bm\sigma-\bm\sigma_h||_{L^2(\Omega\times T)^2}, |||(v,{\bm \sigma}) - ( v_h,{\bm \sigma}_{h})|||_{\text{DG}} \lesssim C(v,\bm\sigma)\mathcal{O}(\rho^{\frac{1}{4}})$, which coincides with the optimal convergence rates of the errors with respect to $\rho$ indicated by the Theorem \ref{anisoerr} for the homogeneous case.

\subsection{Comparisons of TDG and DG methods for 2D space case with nonhomogeneous sources}

We choose the space-time domain $Q=(0,1)^2\times(0,1)$. Consider the exact smooth solution
 \beq
u({\bf x},t) & = & \text{sin}(\pi x_1) \text{sin}(\pi x_2)\text{sin}(\sqrt{3}\pi t),
\cr v& =& \frac{\partial U}{\partial t}, ~ \bm\sigma = -A^{\frac{1}{2}}\nabla U,
\eq
to the IBVP (\ref{model}), with nonhomogeneous soure $f$. We use uniform mesh with $ h_{{\bf x}} \approx h_t =2^{-l}, l\in \mathbb{N}$.

\subsubsection{The case of Dirichlet boundary conditions}

We would like to compare the errors of the approximations generated by the proposed combined DG method of Section 7 and the high-order DG method of \cite{BMPS}. We choose $\Gamma_D=\partial\Omega$. The convergence rates are given in the Table \ref{nnt2dcomdri} for different choices of $p=q+1$. 


\begin{center}
       \tabcaption{}
\label{nnt2dcomdri} \vskip -0.3in
       Convergence rates of the combined space-time DG scheme w.r.t. $h$.
\begin{tabular}{|c|c|c|cc|cc|cc|} \hline
   &   & &  \multicolumn{2}{|c|}{$v_h$} &  \multicolumn{2}{|c|}{$\bm\sigma_h$} &  \multicolumn{2}{|c|}{$|||\cdot|||_{\text{DG}}$} \\ \hline
   $(p,q)$ & \text{Method} & \text{DOFs} & \text{Error} & \text{Rate} & \text{Error} & \text{Rate} &\text{Error} & \text{Rate} \\ \hline
    \multirow{2}*{\((2,1)\)} & \text{TDG} &7680 & 1.30e-2  & 2.79 &  6.33e-3  & 3.61 & 1.57e-1 &  1.86 \\ \cline{2-9}
   & \text{DG} & 12288 & 2.18e-2  & 2.89 & 4.16e-2   & 2.66 & 3.88e-1 & 1.48  \\ \hline
     \multirow{2}*{\((3,2)\)} & \text{TDG} & 5184 &  2.96e-3  & 3.48 & 3.08e-3   & 3.93 & 3.23e-2  & 2.77   \\ \cline{2-9}
   & \text{DG} & 5184   & 8.50e-3   & 3.41 & 6.65e-3  & 4.00 & 1.46e-1  & 2.52   \\ \hline
    \multirow{2}*{\((4,3)\)} & \text{TDG} & 7560 & 2.06e-4 & 4.64 & 2.86e-4  & 4.99 & 2.45e-3   & 3.70  \\ \cline{2-9}
    & \text{DG} & 12288 & 6.05e-4 & 4.01 & 8.28e-4 & 4.14 & 1.42e-2 & 3.53   \\ \hline
   \end{tabular}
     \end{center}


We can see from Table \ref{nnt2dcomdri} that, the convergence orders of errors generated by the global TDG method in $L^2(\Omega\times\{T\})-$norm are between $\mathcal{O}(h^{p+\frac{1}{2}})$ and $\mathcal{O}(h^{p+1})$.
 Besides, the convergence rates of errors in $|||\cdot|||_{\text{DG}}-$norm are $O(h^{\text{min}\{p+\frac{1}{2},q+\frac{1}{2}\}})$, which support convergence rate optimality on the uniform refined meshes of the global TDG method
  just as for the homogeneous case.
Moreover, the approximations generated by the TDG are more accurate than those generated by the high-order DG method, even if the DOFs of the DG method are significantly higher than that of the TDG method.

 \subsubsection{The case of mixed boundary conditions}

We choose $\Gamma_D=\{x_1=0,1\}\times [0,1]$ and $\Gamma_N=[0,1]\times\{x_2=0,1\} $. The convergence rates are given in the Table \ref{nnt2dcommix} for different choices of $p=q+1$.


\begin{center}
       \tabcaption{}
\label{nnt2dcommix} \vskip -0.3in
       Convergence rates of the combined space-time DG scheme w.r.t. $h$.
\begin{tabular}{|c|c|c|cc|cc|cc|} \hline
   &   & &  \multicolumn{2}{|c|}{$v_h$} &  \multicolumn{2}{|c|}{$\bm\sigma_h$} &  \multicolumn{2}{|c|}{$|||\cdot|||_{\text{DG}}$} \\ \hline
   $(p,q)$ & \text{Method} & \text{DOFs} & \text{Error} & \text{Rate} & \text{Error} & \text{Rate} &\text{Error} & \text{Rate} \\ \hline
    \multirow{2}*{\((2,1)\)} & \text{TDG} & 7680  & 1.47e-2 &  3.45  & 8.04e-3   & 2.98 & 1.60e-1 & 1.86  \\ \cline{2-9}
   & \text{DG} &-- &--   &--  & --   &-- &--  & --  \\ \hline
     \multirow{2}*{\((3,2)\)} & \text{TDG} & 5184 & 2.81e-3  & 3.55 & 2.21e-3   & 3.74 & 2.13e-2 &  2.65  \\ \cline{2-9}
   & \text{DG} &-- &--   &--  & --   &-- &--  & --    \\ \hline
    \multirow{2}*{\((4,3)\)} & \text{TDG} & 7560  &  6.66e-4 & 4.53  & 4.94e-4   &  4.92 & 8.05e-3 &  3.56  \\ \cline{2-9}
   & \text{DG} & -- &--   &--  & --   &-- &--  & --  \\ \hline
   \end{tabular}
     \end{center}

 It can be seen from Table \ref{nnt2dcommix} that, the TDG scheme still works, but high-order DG scheme fails (here ``--" represents that the corresponding numerical method has no accuracy) for the case of nonhomogeneous Neumann boundary conditions, which verifies the existing theoretical and numerical conclusions, see Remarks \ref{remarkcomboundhomo}, \ref{remarkcombound} and \cite[Sections 3 and 6]{BMPS}.

\subsubsection{The case of Neumann boundary conditions}

We choose $\Gamma_N=\partial\Omega$. The convergence rates are given in the Table \ref{nnt2dcomneu} for different choices of $p=q+1$. 


\begin{center}
       \tabcaption{}
\label{nnt2dcomneu} \vskip -0.3in
       Convergence rates of the combined space-time DG scheme w.r.t. $h$.
\begin{tabular}{|c|c|c|cc|cc|cc|} \hline
   &   & &  \multicolumn{2}{|c|}{$v_h$} &  \multicolumn{2}{|c|}{$\bm\sigma_h$} &  \multicolumn{2}{|c|}{$|||\cdot|||_{\text{DG}}$} \\ \hline
   $(p,q)$ & \text{Method} & \text{DOFs} & \text{Error} & \text{Rate} & \text{Error} & \text{Rate} &\text{Error} & \text{Rate} \\ \hline
    \multirow{2}*{\((2,1)\)} & \text{TDG} & 7680 & 1.38e-2  & 3.32  & 1.23e-2   &  3.13 & 1.56e-1 &  1.87   \\ \cline{2-9}
   & \text{DG} & --&  -- & -- &   -- & --& -- & --  \\ \hline
     \multirow{2}*{\((3,2)\)} & \text{TDG} & 5184 &  1.43e-3 & 3.95  & 2.11e-3   &  3.42 & 2.22e-2 & 2.53  \\ \cline{2-9}
   & \text{DG} & --&  -- & -- &   -- & --& -- & --    \\ \hline
    \multirow{2}*{\((4,3)\)} & \text{TDG} & 7560 & 5.55e-5  & 4.80  &  6.24e-5  & 5.07  & 8.49e-4 & 3.67   \\ \cline{2-9}
   & \text{DG} &--&  -- & -- &   -- & --& -- & --    \\ \hline
   \end{tabular}
     \end{center}

The conclusion coincides with the above Section.

     \subsection{Discussion on the size of local nonhomogeneous problems} \label{overlapsec}
Consider the exact solution of the Section \ref{test2d1t}. The convergence rates are given in the Table \ref{nnt2dsize} for different space size of local nonhomogeneous problems and different choices of $p=q+1$.


\begin{center}
       \tabcaption{}
\label{nnt2dsize} \vskip -0.3in
       Convergence rates of the combined space-time DG scheme w.r.t. $h$.
\begin{tabular}{|c|c|cc|cc|cc|} \hline
   &   &  \multicolumn{2}{|c|}{$v_h$} &  \multicolumn{2}{|c|}{$\bm\sigma_h$} &  \multicolumn{2}{|c|}{$|||\cdot|||_{\text{DG}}$} \\ \hline
  ($p,q$) & $(\text{size}_{x_1},\text{size}_{x_2})$ & \text{Error} & \text{Rate} & \text{Error} & \text{Rate} &\text{Error} & \text{Rate} \\ \hline
    \multirow{2}*{\((2,1)\)} & (1 1) &1.38e-2  & 3.32  & 1.23e-2   &  3.13 & 1.56e-1 &  1.87 \\ \cline{2-8}
   & $(h^{-1}_{{\bf x}},h^{-1}_{{\bf x}})$ &1.45e-2 & 3.29 &1.09e-2    & 3.28  & 1.67e-1 & 1.95  \\ \hline
        \multirow{2}*{\((3,2)\)} & (1 1) &  1.43e-3 & 3.95  & 2.11e-3   &  3.42 & 2.22e-2 & 2.53  \\ \cline{2-8}
   & $(h^{-1}_{{\bf x}},h^{-1}_{{\bf x}})$  & 3.12e-3  &  3.23  & 3.00e-3   & 3.15  & 3.45e-2 &  2.53  \\ \hline
     \multirow{2}*{\((4,3)\)} & (1 1) & 5.55e-5  & 4.80  &  6.24e-5  & 5.07  & 8.49e-4 & 3.67 \\ \cline{2-8}
   & $(h^{-1}_{{\bf x}},h^{-1}_{{\bf x}})$ & 4.77e-5  & 5.02  &  5.97e-5  & 5.13  & 7.43e-4 & 3.86  \\ \hline
   \end{tabular}
     \end{center}

Here $(\text{size}_{x_1},\text{size}_{x_2})=(1,1)$ represents the involved number of spacial elements employed by the overlapping local DG variational formulation (\ref{localnonhomomacro}) in each coordinate direction. For the case of $(\text{size}_{x_1},\text{size}_{x_2})=(h^{-1}_{{\bf x}},h^{-1}_{{\bf x}})$, it represents that 
the nonoverlapping variant of local DG variational formulation (\ref{localnonhomomacro}) are defined on each nonoverlapping time slab $D_n (1\leq n \leq N)$.
 It can be seen that, if nonhomogeneous local problems are defined on $D_n$, the combined DG scheme can also preserve the same orders of convergence as for the case of $(\text{size}_{x_1},\text{size}_{x_2})=(1,1)$, i.e. the local smooth overlapping case (\ref{localnonhomomacro}). 

\subsection{Heterogeneous media case} \label{heterosec}
Consider the space-time domain $Q=(0,1)^2\times(0,1)$. The anisotropic matrix is the piecewise constant matrix
$$A = \left(
  \begin{array}{cc}
  \frac{3}{4} & \frac{1}{4}  \\
  \frac{1}{4}  & \frac{3}{4}\\
  \end{array}
\right) ~\text{when} ~x_1\leq 0.25,
\quad A = \left(
  \begin{array}{cc}
  \frac{5}{8} & \frac{3}{8}  \\
\frac{3}{8}  & \frac{5}{8} \\
  \end{array}
\right) ~\text{when}~ x_1 > 0.25. $$
As the initial condition, we take a Gaussian wave (see \cite[Section 8.3]{BMPS}) given by $$U_0 = \text{exp}(-||{\bf x}-{\bf x}_0||^2/\zeta^2), \quad v_0=0, \quad {\bm \sigma}_0=-A^{\frac{1}{2}}\nabla U_0,$$
where ${\bf x}_0=(0.5,0.5)^T$ and $\zeta=0.01$. We consider
homogeneous Dirichlet boundary conditions.

\begin{figure}[H]
\begin{center}
\begin{tabular}{ccc}
  \epsfxsize=0.3\textwidth\epsffile{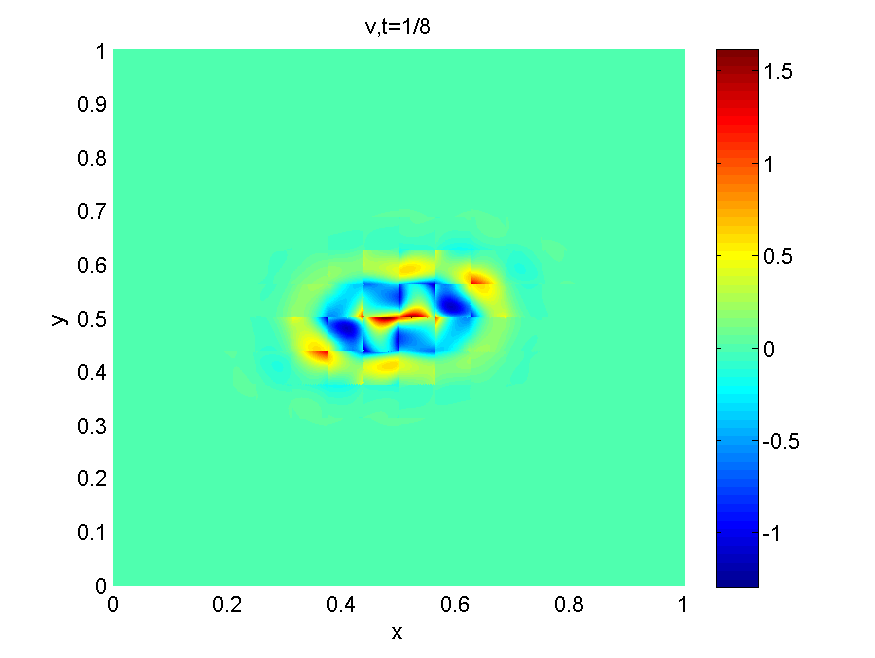}&
 \epsfxsize=0.3\textwidth\epsffile{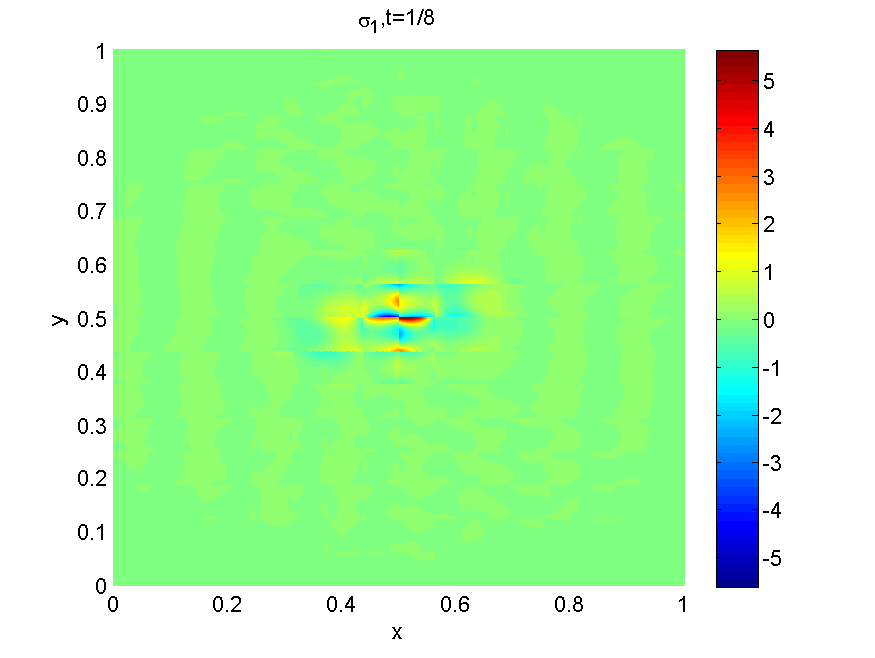}&
  \epsfxsize=0.3\textwidth\epsffile{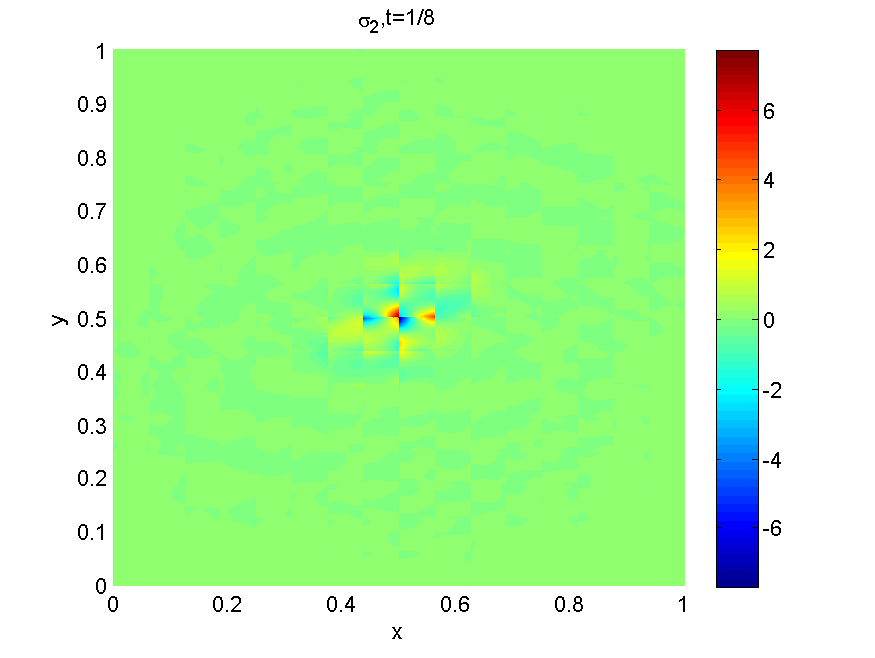}\\
    \epsfxsize=0.3\textwidth\epsffile{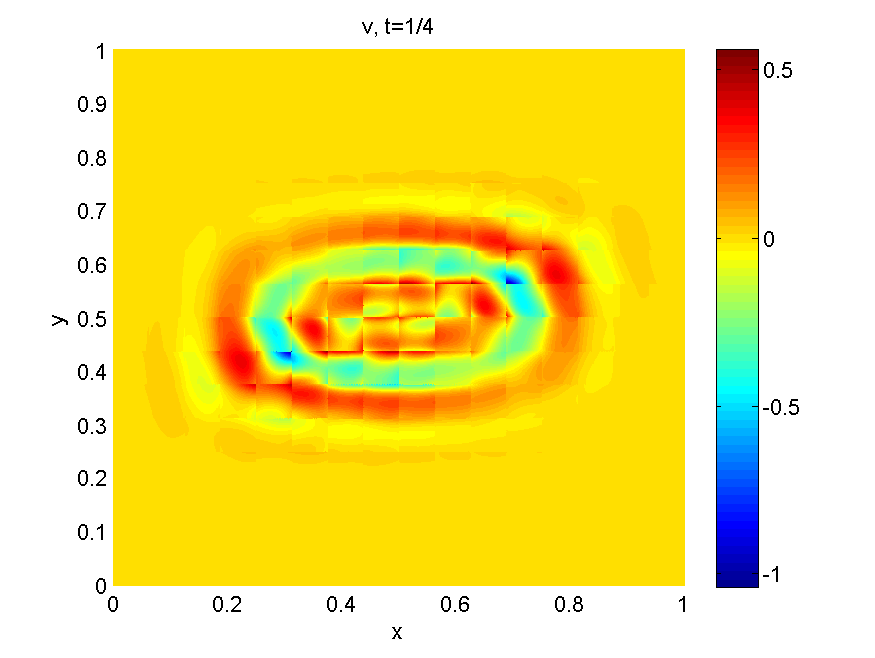}&
 \epsfxsize=0.3\textwidth\epsffile{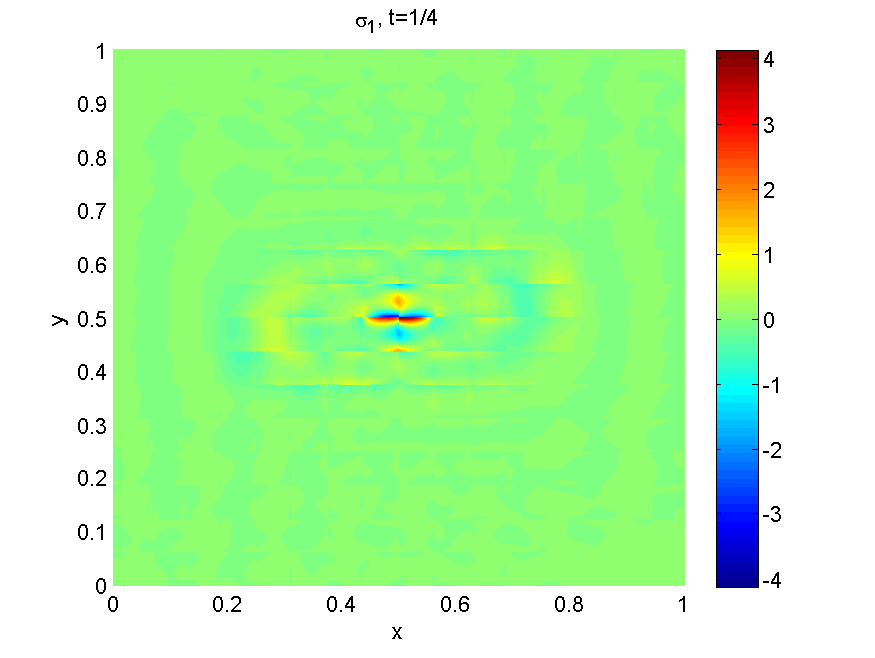}&
  \epsfxsize=0.3\textwidth\epsffile{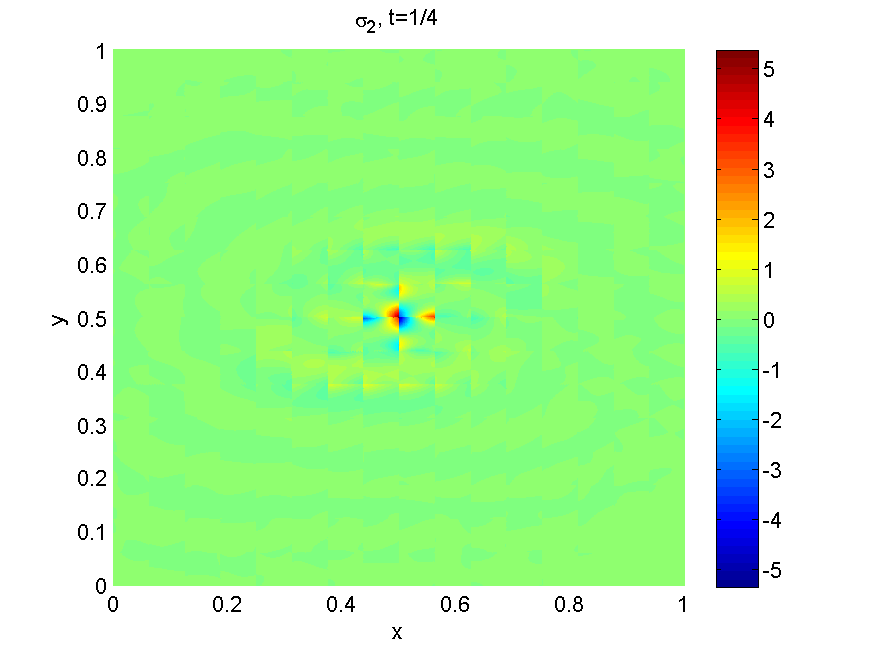}\\
   \epsfxsize=0.3\textwidth\epsffile{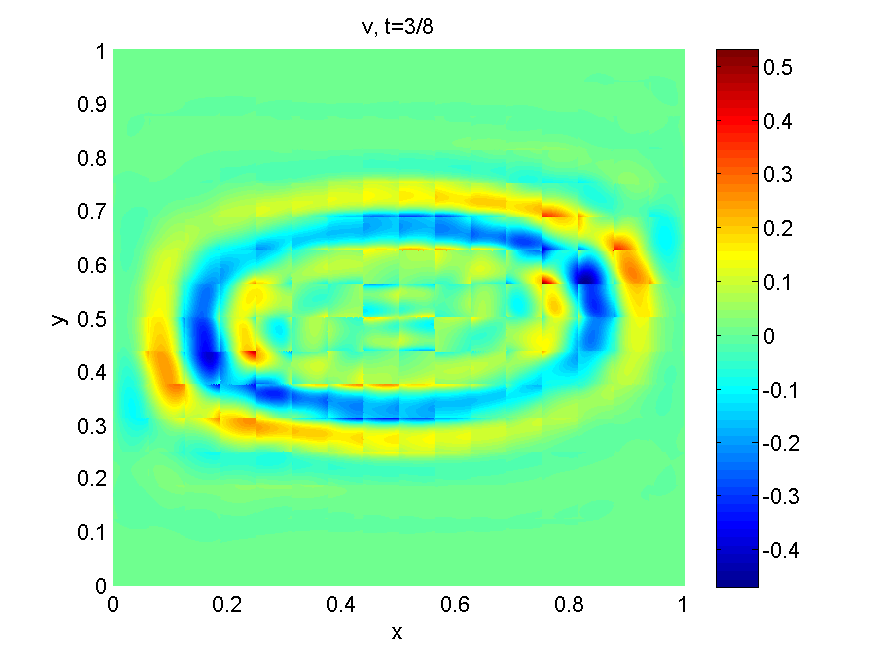}&
 \epsfxsize=0.3\textwidth\epsffile{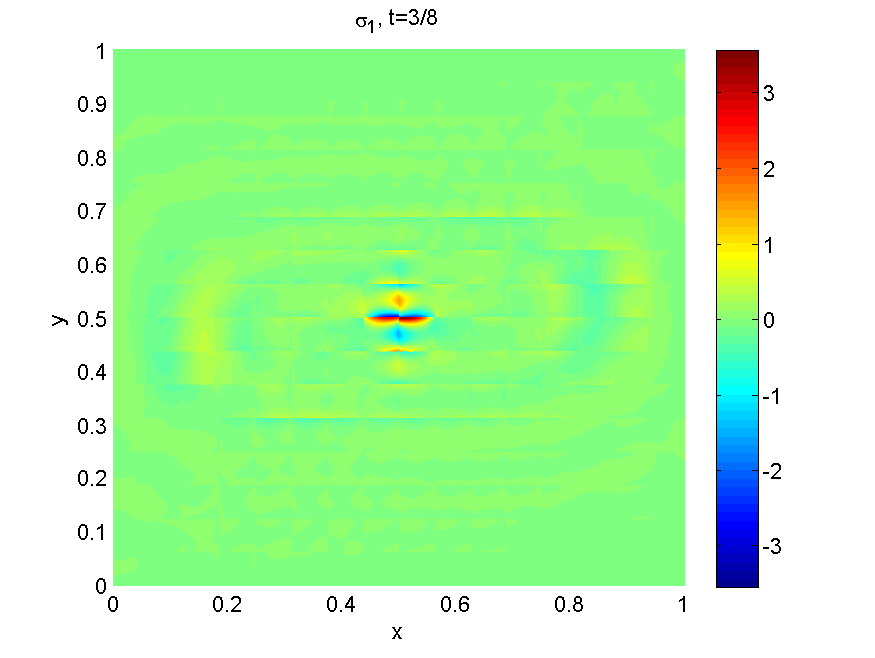}&
  \epsfxsize=0.3\textwidth\epsffile{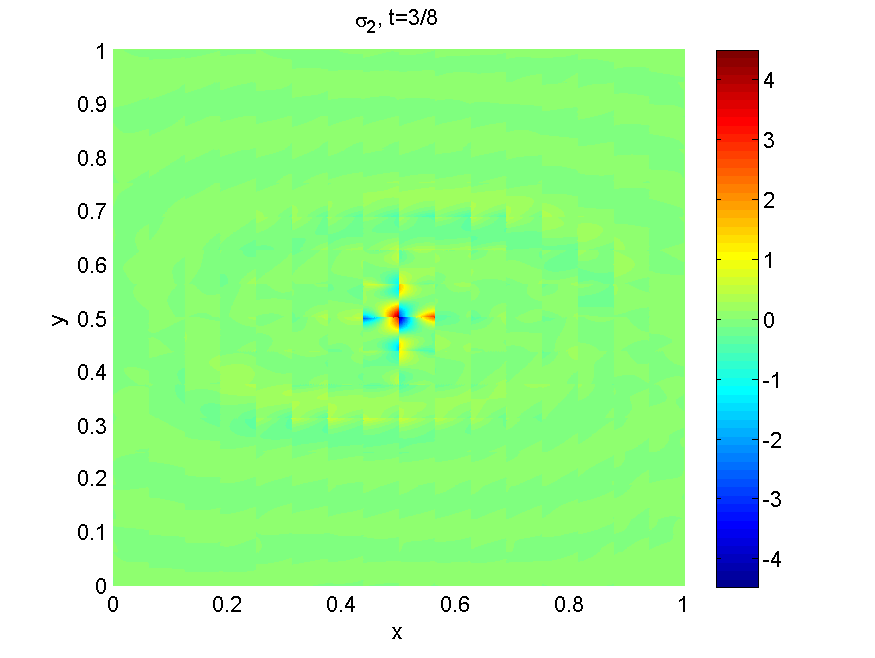}\\
\end{tabular}
\end{center}
  \caption{numerical solution } \label{numericalsolu2}
\end{figure}

Snapshots of the solution are shown in Figure \ref{numericalsolu2}. First, the initial condition evolves in the right homogeneous medium. At time $t = 0.25$, the wave crosses over the interface with the discontinuous anisotropic matrices, and  into the medium with higher wave speed. The snapshot at $t = 0.375$ shows that the incident wave is transmitted across the interface with higher wave speed and a shallow wavefront.

For the case in layered media (i.e., $A$ is not constant), it is difficult
to give an analytic solution of the homogeneous acoustic wave system (\ref{model}). In order to compute accuracies of the Trefftz DG approximations generated by the proposed method, as usual we replace the analytic solution by a good approximation
generated by the same method with very fine grids. The convergence rates are given in the Table \ref{hetennt1} for the case of $p=3$. %

\begin{center}
       \tabcaption{}
\label{hetennt1} \vskip -0.3in
       Convergence rates of the space-time TDG scheme w.r.t. $h$.
\begin{tabular}{|c|c|c|c|c|c|c|c|} \hline
      &  \multicolumn{2}{|c|}{$v_h$} &  \multicolumn{2}{|c|}{$\sigma_h$} &  \multicolumn{2}{|c|}{$|||\cdot|||_{\text{DG}}$} \\ \hline
 $h$ & \text{Error} & \text{Rate} & \text{Error} & \text{Rate} &\text{Error} & \text{Rate} \\ \hline
  1/4 & 1.76e-2  &  & 5.16e-2  &    & 9.26e-2  &      \\ \hline
  1/8 & 1.02e-3 & 4.11  & 2.77e-3  & 4.22   & 8.36e-3  &  3.47   \\ \hline
   1/16 & 6.29e-5  & 4.02 &  1.68e-4 &  4.04  & 7.34e-4  &   3.51   \\ \hline
    \end{tabular}
     \end{center}


From the Table \ref{hetennt1}, we can obtain that $|| v-v_h||_{L^2(\Omega\times T)} \approx ||\bm\sigma-\bm\sigma_h||_{L^2(\Omega\times T)^2} =\mathcal{O}(h^{p+1})$. The last column showing the experimental convergence rates of the errors measured in $|||\cdot|||_{\text{DG}}-$norm indicates that the estimates of Theorem \ref{piecewise_anisoerr} are sharp.

\section{Appendix: the derivation of Eq. (\ref{initanisotoiso})}
For convenience, we use ${\bf p}_1,{\bf p}_2, \cdots, {\bf p}_d$ to denote the column vectors of $P$, and use ${\bf q}_1^T,{\bf q}_2^T,\cdots, {\bf q}_d^T$ to denote the row vectors of $P$. Then each of these vectors is a unit vector, and ${\bf p}_1,{\bf p}_2, \cdots, {\bf p}_d$ (and ${\bf q}_1,{\bf q}_2,\cdots, {\bf q}_d$) are orthogonal each other.

For $1\leq i \leq d$, by the coordinate transformation (\ref{tran1}), we have
\vskip -0.1in
\be \nonumber
\frac{\partial \hat v}{\partial x_i} = \sum_{k=1}^d \frac{\partial \hat v}{\partial \hat x_k} \frac{\partial \hat x_k} {\partial x_i}  = \sum_{k=1}^d \frac{\partial \hat v}{\partial \hat x_k} \lambda_k^{-\frac{1}{2}}q_{ik},
\en
which yields
\be \nonumber 
\nabla \hat v = P^T \Lambda^{-\frac{1}{2}} \hat\nabla \hat v.
\en
Combining it with (\ref{firstorderrela}), yields
\be
A^{\frac{1}{2}}\nabla v = A^{\frac{1}{2}}\nabla \hat v= A^{\frac{1}{2}}P^T \Lambda^{-\frac{1}{2}} \hat\nabla \hat v=  P^T \hat\nabla\hat v.
\en

Next, by (\ref{firstorderrela}), it holds that
\vskip -0.1in
\be \nonumber
A^{\frac{1}{2}}{\bm \sigma} = P^T \Lambda^{\frac{1}{2}} \hat{\bm \sigma} = \sum_{i=1}^d \lambda_i^{\frac{1}{2}}{\bf q}_{i} \hat{\sigma}_i.
\en
By the chain rule, we obtain
\vskip -0.3in
\beq \label{appenfinal}
\nabla\cdot \bigg(A^{\frac{1}{2}}{\bm \sigma}\bigg) &=& \nabla\cdot \bigg( \sum_{i=1}^d \lambda_i^{\frac{1}{2}}{\bf q}_{i} \hat{\sigma}_i \bigg)  = \sum_{j=1}^d \sum_{i=1}^d \lambda_i^{\frac{1}{2}}q_{ji} \frac{\partial \hat{\sigma}_i}{\partial x_j}
    \cr & = &   \sum_{j=1}^d \sum_{i=1}^d \lambda_i^{\frac{1}{2}}q_{ji} \bigg(\sum_{k=1}^d \frac{\partial \hat{\sigma}_i}{\partial \hat x_k} \frac{{\partial \hat x_k}}{{\partial  x_j}} \bigg) = \sum_{j=1}^d \sum_{i=1}^d \lambda_i^{\frac{1}{2}}q_{ji} \bigg(\sum_{k=1}^d \frac{\partial \hat{\sigma}_i}{\partial \hat x_k} \lambda_k^{-\frac{1}{2}} q_{jk} \bigg)
   \cr & = &   \sum_{i=1}^d \sum_{k=1}^d \lambda_i^{\frac{1}{2}}\lambda_k^{-\frac{1}{2}}  \frac{\partial \hat{\sigma}_i}{\partial \hat x_k} \bigg(\sum_{j=1}^d q_{ji}q_{jk} \bigg) = \sum_{i=k=1}^d \frac{\partial \hat{\sigma}_i}{\partial \hat x_k} = \hat\nabla \cdot \hat{\bm \sigma}.
 \eq

\end{document}